\documentclass[a4paper]{article}
\usepackage{amsthm,amsmath,amssymb}
\usepackage{color}

\newtheorem{teor}{Theorem}[section]
\newtheorem{defin}[teor]{Definition}
\newtheorem{lemm}[teor]{Lemma}
\newtheorem{osse}[teor]{Remark}
\newtheorem{prop}[teor]{Proposition}
\newtheorem{defi}[teor]{Definition}
\newtheorem{coro}[teor]{Corollary}
\newtheorem{prob}[teor]{Problem}

\newcommand{\bele}{\begin{lemm}\begin{sl}}
\newcommand{\enle}{\end{sl}\end{lemm}}
\newcommand{\bedef}{\begin{defi}\begin{sl}}
\newcommand{\eddef}{\end{sl}\end{defi}}
\newcommand{\bete}{\begin{teor}\begin{sl}}
\newcommand{\ente}{\end{sl}\end{teor}}
\newcommand{\beos}{\begin{osse}\begin{rm}}
\newcommand{\eddos}{\end{rm}\end{osse}}
\newcommand{\bepr}{\begin{prop}\begin{sl}}
\newcommand{\empr}{\end{sl}\end{prop}}
\newcommand{\bepro}{\begin{prob}\begin{rm}}
\newcommand{\empro}{\end{rm}\end{prob}}
\newcommand{\bede}{\begin{defin}\begin{sl}}
\newcommand{\edde}{\end{sl}\end{defin}}
\newcommand{\beco}{\begin{coro}\begin{sl}}
\newcommand{\enco}{\end{sl}\end{coro}}

\newcommand{\quext}{\quad\text}
\newcommand{\qquext}{\qquad\text}

\newcommand{\RR}{\mathbb{R}}

\newcommand{\beeq}[1]{\begin{equation}\label{#1}}
\newcommand{\eddeq}{\end{equation}}

\newcommand{\beeqa}[1]{\begin{eqnarray}\label{#1}}
\newcommand{\eddeqa}{\end{eqnarray}}

\newcommand{\beal}[1]{\begin{align}\label{#1}}
\newcommand{\eddal}{\end{align}}

\newcommand{\bespl}[1]{\begin{split}\label{#1}}
\newcommand{\edspl}{\end{split}}

\newcommand{\bega}[1]{\begin{gather}\label{#1}}
\newcommand{\edga}{\end{gather}}

\newcommand{\beeqax}{\begin{eqnarray*}}
\newcommand{\eddeqax}{\end{eqnarray*}}

\def\qed{\ifmmode % if math mode, assume display: omit penalty etc.
  \else \leavevmode\unskip\penalty9999 \hbox{}\nobreak\hfill
  \fi
  \quad\hbox{\hskip.5em\vrule width.4em height.6em depth.05em\hskip.1em}}
\def\endproofsym{\qed}
\renewenvironment{proof}[1][Proof]{\trivlist\item[\hskip\labelsep{\hskip0pt
    {\normalfont\scshape#1.}\hskip .321429\parindent}]\ignorespaces}
{\endproofsym\endtrivlist}
\def\endnobox{\def\endproofsym{}\end{proof}\def\endproofsym{\qed}}

\newcommand{\no}{\nonumber}

\newcommand{\beeqao}{\begin{eqnarray}\no}
\newcommand{\bealo}{\begin{align}\no}
\newcommand{\besplo}{\begin{split}\no}
\newcommand{\begao}{\begin{gather}\no}

\newcommand{\duav}[1]{\langle{#1}\rangle}

\newcommand{\io}{\int_\Omega}

\newcommand{\epsi}{\varepsilon}
\newcommand{\ee}{_{\varepsilon}}

\newcommand{\OO}{_{\Omega}}

\newcommand{\fhi}{\varphi}

\newcommand{\lhs}{left hand side}
\newcommand{\rhs}{right hand side}

\DeclareMathOperator{\per}{per}
\DeclareMathOperator{\dive}{div}
\DeclareMathOperator{\deriv}{d}

\DeclareMathOperator{\sign}{sign}

\newcommand{\HUH}{H^1(0,T;H)}

\newcommand{\HUVp}{H^1(0,T;V')}

\newcommand{\LDH}{L^2(0,T;H)}
\newcommand{\LDV}{L^2(0,T;V)}
\newcommand{\LDVp}{L^2(0,T;V')}
\newcommand{\LIH}{L^\infty(0,T;H)}
\newcommand{\LIV}{L^\infty(0,T;V)}

\newcommand{\LIVp}{L^\infty(0,T;V')}

\newcommand{\LDHD}{L^2(0,T;H^2(\Omega))}

\let\TeXchi\chi
\def\chi{{\setbox0 \hbox{\mathsurround0pt
$\TeXchi$}\hbox{\raise\dp0 \copy0 }}}

\newcommand{\teta}{\vartheta}

\newcommand{\calI}{{\mathcal I}}

\newcommand{\calT}{{\mathcal T}}

\newcommand{\calE}{{\mathcal E}}

\newcommand{\calN}{{\mathcal N}}

\newcommand{\calB}{{\mathcal B}}
\newcommand{\calO}{{\mathcal O}}

\newcommand{\ubarvt}{\underline{\vt}}

\newcommand{\baru}{\overline{u}}

\newcommand{\dit}{\deriv\!t}
\newcommand{\dis}{\deriv\!s}

\newcommand{\ddt}{\frac{\deriv\!{}}{\dit}}

\newcommand{\zee}{_{0,\varepsilon}}

\newcommand{\vt}{\vartheta}

\newcommand{\vu}{\boldsymbol{u}}

\newcommand{\barvt}{\overline{\vt}}

\newcommand{\barvu}{\overline{\vu}}

\newenvironment{bettirev}{\color{blue}}{\color{black}}
\newcommand{\bber}{\begin{bettirev}}
\newcommand{\eber}{\end{bettirev}}

\newenvironment{michelarev}{\color{red}}{\color{black}}
\newcommand{\bmicr}{\begin{michelarev}}
\newcommand{\emicr}{\end{michelarev}}

\numberwithin{equation}{section}
\begin{document}

\title{Existence of solutions to a two-dimensional 
  model for nonisothermal 
  two-phase flows of incompressible fluids}

\author{
{Michela Eleuteri\thanks{Dipartimento di Matematica ed Informatica ``U. Dini", viale Morgagni 67/a, 
I-50134 Firenze, Italy.
E-mail: \textit{eleuteri@math.unifi.it}. 
The author is partially supported by the FP7-IDEAS-ERC-StG 
   Grant~\#256872 (EntroPhase) and by GNAMPA (Gruppo Nazionale per l'Analisi Matematica, 
   la Probabilit\`a e le loro Applicazioni) of INdAM (Istituto Nazionale di Alta Matematica)}}
\and
{Elisabetta Rocca\thanks{Weierstrass Institute for Applied
    Analysis and Stochastics, 
 Mohrenstrasse~39, D-10117 Berlin, Germany. 
E-mail: \textit{rocca@wias-berlin.de} and Dipartimento 
   di Matematica ``F. Enriques'',
Universit\`a degli Studi di Milano, I-20133 Milano, Italy.
E-mail: \textit{elisabetta.rocca@unimi.it}. 
The author is supported by the FP7-IDEAS-ERC-StG
   Grant~\#256872 (EntroPhase)  and by GNAMPA (Gruppo Nazionale per l'Analisi Matematica, la Probabilit\`a e le loro Applicazioni) of INdAM (Istituto Nazionale di Alta Matematica)}}
\and
{Giulio Schimperna\thanks{Dipartimento di Matematica ``F. Casorati'', 
Universit\`a degli Studi di Pavia,
via Ferrata 1, I-27100 Pavia, Italy.
E-mail: \textit{giusch04@unipv.it}. 
The author is partially supported by the FP7-IDEAS-ERC-StG
   Grant~\#256872 (EntroPhase) and by GNAMPA (Gruppo Nazionale per l'Analisi Matematica, la Probabilit\`a e le loro Applicazioni) of INdAM (Istituto Nazionale di Alta Matematica)}}}

\maketitle

\begin{abstract}\noindent
We consider a thermodynamically consistent diffuse interface model 
describing two-phase flows of incompressible fluids in a non-isothermal
setting. The model was recently introduced in~\cite{ERS} where existence
of weak solutions was proved in three space dimensions. Here, we 
aim at studying the properties of solutions in the two-dimensional
case. In particular, we can show existence of global in time 
solutions satisfying a stronger formulation of the model with
respect to the one considered in~\cite{ERS}. Moreover, we can admit 
slightly more general conditions on some material coefficients
of the system.
%Moreover,
%we analyze the long-time behavior of the model proving existence
%of the global attractor in a suitable phase space.
%
\end{abstract}

\smallskip


\smallskip


%%%%%%%%%%%%%%%%%%%%%%%%%%%%%%%%%%%%%%%%%%%%%%%%%%%%%%%%%%%%%%%%%%%%%%%%%%%%%%%%%%%%%%%%%%%%%%%%%
%%%%%%%%%%%%%%%%%%%%%%%%%%%%%%%%%%%%%%%%%%%%%%%%%%%%%%%%%%%%%%%%%%%%%%%%%%%%%%%%%%%%%%%%%%%%%%%%%

\section{Introduction}
\label{sec:intro}

We consider here a mathematical model for two-phase flows
of non-isothermal incompressible fluids in a bounded
container $\Omega$ in $\RR^2$ during a finite time interval $(0,T)$. 
The model consists in a PDE system describing the evolution 
of the unknown variables $\vu$ (macroscopic
velocity), $\fhi$ (order parameter), $\mu$ (chemical potential),
$\vt$ (absolute temperature), and it takes the form
\begin{align}\label{incom}
  & \dive \vu = 0, \\
 \label{ns}
  & \vu_t + \vu \cdot \nabla \vu + \nabla p
    = \Delta \vu - \dive ( \nabla \fhi \otimes \nabla \fhi ),\\
 \label{CH1}
  & \fhi_t + \vu \cdot \nabla \fhi = \Delta \mu, \\
 \label{CH2}
  & \mu = - \Delta \fhi + F'(\fhi) - \vt, \\
 \label{calore}
  & \vt_t + \vu \cdot \nabla \vt + \vt \big( \fhi_t + \vu \cdot \nabla \fhi \big) 
    - \dive(\kappa(\vt)\nabla \vt) = | \nabla \vu |^2 + | \nabla \mu |^2.
\end{align}  
%
%Here and below, $\nabla \vu := ( \nabla \vu + (\nabla \vu)^t )/2$ represents
%the symmetrized velocity gradient. 
Relation~\eqref{ns}, with the 
incompressibility constraint \eqref{incom},
represents a variant of the Navier-Stokes system; 
\eqref{CH1}-\eqref{CH2} correspond 
to a form of the Cahn-Hilliard system \cite{CH} for phase separation, 
while \eqref{calore} is the internal energy equation describing
the evolution of temperature. Note that transport effects are 
admitted for all variables in view of the occurrence of 
material derivatives in \eqref{ns}, \eqref{CH1}, 
and~\eqref{calore}. As usual, the variable $p$ 
in the Navier-Stokes system~\eqref{ns} 
represents the (unknown) pressure. 
The function $F$ whose derivative
appears in \eqref{CH2} is a possibly non-convex
potential whose minima represent the least energy 
configurations of the phase variable. Here we will
assume that $F$ is smooth and has at most a power-like growth
at infinity. Indeed, it is not clear whether our result
could be extended to other classes of physically significant
potential, having nonsmooth or singular character
(like, e.g., the so-called logarithmic potential
$F(r) = (1+r)\log(1+r) + (1-r)\log(1-r) - r^2$ which 
typically appears in Cahn-Hilliard-based models,
see, e.g., \cite{MZ}). Finally, the coefficient
$\kappa(\vt)$ in~\eqref{calore} stands for the heat conductivity
of the fluid. Here we shall assume that 
$\kappa$ grows at infinity like a sufficiently high 
power of $\vt$ (see \eqref{hp:kappa} below).

The PDE system \eqref{incom}-\eqref{CH2} in the case of a constant temperature $\teta$, refereed 
in the literature as Model~H, is a diffuse interface model for incompressible isothermal 
two-phase flows which consists of the Navier-Stokes equations for the (averaged) velocity 
$\vu$ nonlinearly coupled with a convective Cahn-Hilliard equation for the phase variable $\varphi$
(cf., for instance, \cite{AMW,GPV,HMR,HH,JV,Kim2012,LMM} and the references 
\cite{A1,A2,B,CG,GG1,GG2,GG3,LS,S,ZWH,ZF} for the study of the resulting evolution system). 
However, even if many authors considered the isothermal Model H, up to our knowledge no 
contributions are so far  present in the literature in the non-isothermal case, except 
for~\cite{SunLiuXu} where a linearization of \eqref{calore} is considered and our 
previous paper~\cite{ERS} where we introduced system \eqref{incom}-\eqref{calore} in the 3D case. 

The above model was indeed proposed in our recent work~\cite{ERS} starting
from the balance laws for internal energy and entropy; 
in particular, thermodynamical consistence was shown to hold
for any (positive) value of the absolute temperature $\vt$. 
Moreover, existence of solutions for a {\sl weak formulation}\/ 
of \eqref{incom}-\eqref{calore} was proved 
when the system is settled in a smooth bounded domain 
$\Omega\subset \RR^3$ and complemented with no-flux conditions
for $\fhi$, $\mu$ and $\vt$ and with {\sl slip}\/ conditions
for $\vu$. 
%The argument given in \cite{ERS} is based on a 
%suitable reformulation of the model along the lines
%of an idea originally developed in~\cite{BFM09, FEI09}
%for dealing with heat conduction phenomena in fluids 
%and in \cite{FPR09} for solid-liquid phase transitions.
Mathematically speaking, the main source of difficulty 
in system~\eqref{incom}-\eqref{calore}
comes from the quadratic terms on the \rhs\ of 
\eqref{calore}. Their occurrence is physically 
motivated as one considers the derivation of the model in terms
of the energy and entropy balances (cf.~\cite[Sec.~2]{ERS}). 
Roughly speaking, one can say that these terms represent
a source of thermal energy coming from the dissipation of 
kinetic energy due to viscosity (cf.~\eqref{ns}) and 
of configuration energy due to action of micro-forces 
(cf.~\eqref{CH1}-\eqref{CH2}). This energy dissipation, as 
expected, happens in such a way to increase the entropy of the
system.
 
From the analytical viewpoint, the quadratic terms in \eqref{calore}
can be controlled only in the $L^1$-norm, at least in the 3D-case. 
For this reason, proving existence for the formulation 
\eqref{incom}-\eqref{calore} appears to be out of reach.
Actually, the notion of weak solution considered in~\cite{ERS}
is based on a suitable reformulation of the model along the lines
of an idea originally developed in~\cite{BFM09, FEI09}
for dealing with heat conduction in fluids,
in \cite{FPR09} for solid-liquid phase transitions,
and more recently in~\cite{RoRo} for damage phenomena.
In such a setting, the ``heat'' equation \eqref{calore} is 
replaced with a relation describing the balance of 
{\sl total energy}\/ (i.e., not only of thermal energy),
which does no longer contain quadratic terms.
This is complemented with a distributional version of 
the {\sl entropy inequality}. It is worth observing that
the weak formulation considered in \cite{ERS}
is consistent with the standard (strong) 
one \eqref{incom}-\eqref{calore}.
Indeed, it is not difficult to prove that, 
at least for {\it sufficiently smooth}\/ weak solutions, 
the total energy balance together with the entropy
inequality imply the original form of the heat 
(or, more precisely, internal energy
balance) equation \eqref{calore}. However, as noted
above, the required regularity in the 3D case is not at all known.

Looking at the 2D model, whose analysis is the aim of this paper,
it is well-known that, for the Navier-Stokes system \eqref{ns}, 
additional regularity is available provided that the 
forcing term (here given by $-\dive(\nabla\fhi\otimes\nabla\fhi)$)
lies in~$L^2$ (cf., e.g., \cite{rob}). 
Fortunately, this seems to happen in our case,
as one can readily check starting from the available
energy and entropy estimates; hence, there is hope to get
additional summability for the quadratic term $|\nabla \vu|^2$ 
in~\eqref{calore}. This was the motivation which led us 
to investigate whether it is possible to prove existence of 
a solution to the {\sl original}\/ (strong) system
\eqref{incom}-\eqref{calore} 
in two space dimensions. Indeed, we may give 
a positive answer to this question, but the argument 
we use for arriving at this conclusion
is far from being a straighforward one. So, let us
try to give some ideas of the mathematical difficulties we met.

To make things clear, we start by 
introducing some basic assumptions.
First of all, in order to avoid technical 
complications related with the choice of boundary data, 
we ask the system to be settled in the unit torus 
$\Omega:=[0,1]\times [0,1]$ and complemented 
with periodic boundary conditions for all unknowns. 
It is worth noting that, at the price of some notational change
and of limited technical complications, 
other types of boundary conditions could be 
assumed. For instance, we may take no-flux 
(i.e., homogeneous Neumann)
conditions for $\fhi$, $\mu$ and $\vt$ (as it is
physically reasonable if one assumes the container $\Omega$
to be insulated from the exterior), whereas for $\vu$ 
we may consider any conditions that could 
allow the transport terms to have zero spatial 
mean and are compatible with the existence of 
smooth solutions to the 2D Navier-Stokes system
(cf.~estimate~\eqref{st51} below). This is the 
case, for instance, of homogeneous Dirichlet conditions
(cf., e.g., \cite[Thm.~3.10, p.~314]{Tem}).
The evolution is assumed to take place on a given
reference interval $(0,T)$, with no restrictions on
the magnitude of the final time $T>0$. 
%It is worth noting 
%that, although different conditions were considered 
%in~\cite{ERS} for the three-dimensional
%problem, it was also noticed there that the periodic case 
%could be treated with the same methods as well. 

Coming to our mathematical argument, once additional regularity 
for $\vu$ has been obtained, we need to get further bounds
for the other variables, with the aim of proving an estimate
for the remaining nonlinear terms in the system 
(and, particularly, for the quadratic term $| \nabla \mu |^2$
on the \rhs\ of \eqref{calore}). Actually,
due to the strong coupling between the energy balance 
equation \eqref{calore} and the Cahn-Hilliard
system~\eqref{CH1}-\eqref{CH2}, getting a regularity estimate
for the variables $\fhi$, $\mu$ and $\vt$ requires to 
manage all equations simultaneously in a non-straighforward
way. This further regularity estimate represents,
in our view, the main novelty of the present paper. 

Referring to Section~\ref{sec:well} for more details, 
we just give here some brief explanation of the procedure.
The first thing one can naturally do is to differentiate 
in time the Cahn-Hilliard relation \eqref{CH2}, 
and to test the result by $\varphi_t$. However, this trick 
works only if one is able to control the product 
$\teta_t\varphi_t$. In view of the highly nonlinear
structure of \eqref{calore}, getting a bound for $\teta_t$ 
by working directly on the ``heat'' equation seems
difficult. Hence, the only possibility seems 
that of testing \eqref{calore}
by $\varphi_t$ in order to cancel the bad term. 
However, we then need to control the 
quadratic terms $|\nabla \vu|^2$ and $|\nabla\mu|^2$
on the right hand side, and particularly the latter one,
for which only an $L^1$-estimate is available at this level.
Actually, in order to provide a bound for $|\nabla\mu|^2$,
we have to use some duality technique,
and, more precisely, we need to rely on a sharp
two-dimensional interpolation-embedding inequality
which is proved in the Appendix. This property
follows from well-known two-dimensional embedding
theorems; however we were not able to find it anywhere in 
the literature. The underlying idea 
stands in optimizing with respect to $q$ the embedding
constant of the immersion 
$\|v\|_{H^1(\Omega)'}\leq c_q \|v\|_{L^q(\Omega)}$ 
which holds true in 2D for every $q\in (1,+\infty]$ and 
$v\in L^q(\Omega)$ (cf.~\cite[(17), p.~479]{Trud}). 
Then, applying the embedding inequality 
to the term $|\nabla\mu|^2$, and performing a notable
amount of technical work, we can actually prove
the desired enhanced a-priori bounds.
These estimates permit us to pass to the limit 
in a suitable approximation scheme (which is just 
sketched, for brevity), obtaining in this way 
a solution to the original (strong) system
\eqref{incom}-\eqref{calore}, in a proper
regularity class, coupled with periodic boundary 
conditions and with the initial conditions.

It is finally worth noting that, while
in \cite{ERS} we needed to assume non-constant specific 
heat and heat conductivity, both having a 
suitable growth at 0 and at $\infty$ (the reasons
were mainly of mathematical type;
namely, we needed to get a sufficient summability
for $\teta$), here we just need a 
sufficiently fast growing heat conductivity, but 
we can allow for a constant specific heat. 
Even if this choice is mainly motivated by mathematical reasons,
a physical justification for it can be found, e.g., in \cite{Ze}.

The arguments given in this paper may also be adapted 
to deal with other interesting related models. For instance, 
we could consider the case when the ``Cahn-Hilliard'' 
relations \eqref{CH1}-\eqref{CH2}
are replaced by their ``Allen-Cahn'' equivalent 
(cf.~\cite{CA77})
\begin{equation}\label{ac}
  \varphi_t
   + \vu\cdot\nabla\varphi
   - \Delta\varphi
   + F'(\varphi)
   -\teta = 0.
\end{equation}
Moreover, as a byproduct of our results, one can deduce the existence of 
solutions in 2D for the so-called Fr\'emond's model of phase transitions 
with microscopic effects introduced in~\cite{Frem},
at least in the case of power-like heat conductivity.
The Fr\'emond model basically corresponds to system (\ref{incom})-(\ref{calore}) 
where the velocity $\vu$ is assumed to be 
identically equal to 0. Indeed, for this model,
in the case of Neumann boundary conditions and 
standard Fourier heat flux law, existence of global in time 
``strong'' solutions was known
only in the one-dimensional setting (cf.~\cite{ls1,ls2}),
while {\sl weak solutions}\/ were proved to exist in 3D (cf.~\cite{FPR09})
when \eqref{calore} is replaced by the {\sl total energy balance}\/ and an {\sl entropy inequality}.
Hence, this paper covers the missing 2D case, at least for the 
case of power-like heat conductivity and periodic boundary conditions.

Let us finally note that uniqueness of solutions, as well as their long-time behavior 
(both in terms of trajectories and of attractors) for the whole 
system~\eqref{incom}-\eqref{calore},
are still open issues, which will be the subject of further investigations.

Here is the plan of the paper: in the next
Section~\ref{sec:main} we specify our assumptions on coefficients
and data and state the precise mathematical formulation of 
our problem together with the related existence 
theorem. The remainder of the paper is devoted to the proof
of the theorem. In particular, the core of our argument is given 
in Section~\ref{sec:well}, where we provide the a-priori 
estimates and the compactness argument necessary to pass to
the limit in the approximation scheme. Indeed, in order 
to avoid technicalities, the estimates are obtained in a formal
way leaving the details of a possible regularization  
in the subsequent Section~\ref{sec:appro}. Finally, the
Appendix contains the proof of the
mentioned two-dimensional interpolation-embedding
inequality, which plays a key role in the derivation of
the a-priori bounds.

%%%%%%%%%%%%%%%%%%%%%%%%%%%%%%%%%%%%%%%%%%%%%%%%%%%%%%%%%%%%%%%%%%%%%%%%%%%%%%%%%%%%%%%%%%%%%%%%%
%%%%%%%%%%%%%%%%%%%%%%%%%%%%%%%%%%%%%%%%%%%%%%%%%%%%%%%%%%%%%%%%%%%%%%%%%%%%%%%%%%%%%%%%%%%%%%%%%

\section{Assumptions and main results}
\label{sec:main}

In order to state the precise mathematical formulation of 
our problem we first need to introduce some functional spaces.
Recalling that $\Omega=[0,1]\times[0,1]$, 
we note as $H:=L^2_{\per}(\Omega)$ the space of 
functions in $L^2(\RR^2)$ which are $\Omega$-periodic 
(i.e., $1$-periodic both in $x_1$ and in $x_2$). Analogously,
we set $V:=H^1_{\per}(\Omega)$. The spaces $H$ and $V$
are endowed with the norms of $L^2(\Omega)$ 
and $H^1(\Omega)$, respectively. For brevity, the norm
in $H$ will be simply indicated by $\| \cdot \|$. 
We will note by $\| \cdot \|_X$ the norm
in the generic Banach space $X$. The symbol
$\duav{\cdot,\cdot}$ will indicate the duality
between $V'$ and $V$ and $(\cdot,\cdot)$ 
will stand for the scalar product of~$H$. 
We also write $L^p(\Omega)$ in place of $L^p_{\per}(\Omega)$,
and the same for other spaces; indeed, no confusion should
arise since periodic boundary conditions are assumed
to hold for all unknowns. Still for brevity, we use the
same notation for indicating vector-valued (or tensor-valued)
function spaces and related norms. For instance, 
writing $\vu \in H$, we will in fact mean 
$\vu\in L^2_{\per}(\Omega)^2$. Also the incompressibility 
constraint \eqref{incom} will not be emphasized in the 
notation for functional spaces (hence, the notation $\vu \in H$
will also implicitly subsume that $\dive \vu = 0$ in the 
sense of distributions). These simplifications will
allow us to shorten a bit some formulas.

For any function $v\in H$, we will note as
\begin{equation}\label{voo}
  v\OO:= \frac1{|\Omega|} \io v = \io v
\end{equation}
the spatial mean of $v$. Replacing the integral with a duality
pairing, the same notation will be used in case $v\in V'$.
The symbols $V_0$, $H_0$ and $V_0'$ denote the subspaces
of $V$, $H$ and, respectively, $V'$ containing the function(al)s
having zero spatial mean. We notice that the 
distributional operator $(-\Delta)$ is 
invertible if seen as a mapping from $V_0$ to $V_0'$.
In the sequel we shall denote by $\calN$ its inverse 
operator.

Moreover, in the following we will frequently use the following
2D interpolation inequalities:
\begin{align}\label{dis:L4}
  & \| v \|_{L^4(\Omega)} 
   \le c \| v \|_{V}^{1/2} \| v \|^{1/2},\\
 \label{dis:Linfty}
  & \| v \|_{L^\infty(\Omega)} 
   \le c \| v \|_{H^2(\Omega)}^{1/2} \| v \|^{1/2},\\
   \label{dis:Lp}
  & \| v \|_{L^r(\Omega)} 
   \le c \| v \|_{L^s(\Omega)}^{1-\alpha} \| v \|_{L^{\infty}(\Omega)}^{\alpha}, 
   \qquad \alpha = 1 - \frac{s}{r},
\end{align}  
holding for any sufficiently smooth function $v$ and for
suitable embedding constants, all denoted by
the same symbol $c>0$ for brevity.

We will also use the following nonlinear version of the Poincar\'e
inequality 
\begin{equation}\label{poinc}
  \| v^{p/2} \|_V^2
   \le c_p \big( \| v \|_{L^1(\Omega)}^p
    + \| \nabla v^{p/2} \|^2 \big),
\end{equation}
holding for all nonnegative $v\in L^1(\Omega)$ such that
$\nabla v^{p/2}\in L^2(\Omega)$, and for all $p\in [2, \infty)$.
We also note that
\begin{equation}\label{inter-gen} 
  \| v \| \le c \| \nabla v \|^{1/2} \| v \|_{V'}^{1/2}
   \quext{for all }\, v\in V_0.
\end{equation}
This property can be proved by combining the standard interpolation 
inequality $\|v\| \le c \| v \|_V^{1/2} \| v \|_{V'}^{1/2}$
with the Poincar\'e-Wirtinger inequality. 

In the sequel we will frequently use the 
continuous embedding $V\subset L^p(\Omega)$, holding
for all $p\in[1,\infty)$. Actually, using interpolation of
$L^p$-spaces and Young's inequality, it is not difficult to
see that \eqref{inter-gen} implies
\begin{equation}\label{inter-gen2} 
  \| v \|_{L^p(\Omega)}^2 
    \le \epsilon \| \nabla v \|^2 + c_\epsilon \| v \|_{V'}^2
   \quext{for all }\, v\in V_0,
\end{equation}
for all (small) $\epsilon>0$ and correspondingly large $c_\epsilon>0$
whose value additionally depends on $p\in [1,\infty)$.

\smallskip

With the above notation at disposal, we can present
our main assumptions. First of all, we ask the 
configuration potential~$F$ to satisfy:
\begin{align}\label{hp:F1}
  & F\in C^2(\RR;\RR), \quad 
   \liminf_{|r|\to \infty} \frac{F(r)}{|r|} > 0,\\
 \label{hp:F2}
   & F''(r) \ge - \lambda
    \quext{for some }\,\lambda \ge 0
    \quext{and all }\, r\in \RR,\\
 \label{hp:F3}
   & | F''(r) | \le c_F \big( 1 + |r|^{p_F} \big)
    \quext{for some }\,c_F \ge 0,~p_F \ge 0,
    \quext{and all }\, r\in \RR.
\end{align}  
In other words, we ask for $F$ to be a smooth,
coercive (in view of \eqref{hp:F1}),
$\lambda$-convex (cf.~\eqref{hp:F2}) function, 
with at most polynomial growth at infinity (cf.~\eqref{hp:F3}).
These conditions may probably be
relaxed (admitting, for instance, functions with
exponential growth at infinity), at the price, however,
of technical complications. On the other hand, it is not
clear whether it would be possible to admit 
{\sl singular potentials}\/ like the 
logarithmic function mentioned in the Introduction.
%(we refer the reader to \cite{MZ} for a detailed analysis
%of the Cahn-Hilliard equation with singular function~$F$,
%also covering the 2D-case).

Next, we assume the heat conductivity to be given by
\begin{equation}\label{hp:kappa}
  \kappa(r) = 1 + r^q, \quad
    q \in [2,\infty), \quad
    r \ge 0.
\end{equation}
Correspondingly, we define
\begin{equation}\label{defiK}
  K(r) := \int_0^r \kappa(s)\,\dis
     = r + \frac{1}{q+1} r^{q+1}, \quad
     r \ge 0.
\end{equation}
In the sequel we will often need to estimate the value $\|K(\vt)\|^2_V$. 
To this aim, we first observe that, for some $k_q > 0$,
\begin{equation}\label{K11}
  \io \kappa(\vt)^2 | \nabla \vt |^2
   = \| \nabla K(\vt) \|^2 
   \ge \| \nabla \vt \|^2 + k_q \| \nabla \vt^{q+1} \|^2.
\end{equation}
Then, exploiting \eqref{poinc} with the choice $p = 2$ we obtain
\begin{equation}\label{K12}
  \| K(\vt) \|^2_V 
   \le c_q \left (\int_{\Omega} (\vt + \vt^{q+1}) \right )^2 
     + c_q \left (\int_{\Omega}\left | \nabla \vt + \nabla \left (\frac{\vt^{q+1}}{q+1} \right )\right |^2 \right ) 
   =: I + II
\end{equation}
for some $c_q>0$. Now, using again \eqref{poinc}, this time with the choice $p = 2(q+1)$,
we deduce
\begin{equation}\label{K13}
  I \le c_q \| \vt \|_{L^1(\Omega)}^2 
   + c_q \left( \|\vt\|^{2(q+1)}_{L^1(\Omega)} + \|\nabla \vt^{q+1}\|^2 \right)
  \le c_q  \left( 1 + \|\vt\|^{2(q+1)}_{L^1(\Omega)} + \|\nabla \vt^{q+1}\|^2 \right).
\end{equation}
Estimating $II$ with the help of \eqref{K11} we then conclude
%
%On the other hand
%
%\begin{equation}\label{K11}
%  II \le c \int_{\Omega} \kappa^2(\vt) |\nabla \vt|^2.
%\end{equation}
%
\begin{equation}\label{stima:Kvt}
  \| K(\vt) \|^2_V 
   \le c_q \left (1 + \|\vt\|^{2(q+1)}_{L^1(\Omega)}
        + \int_{\Omega} \kappa^2(\vt) |\nabla \vt|^2 \right ).
\end{equation}
Finally we come to our assumptions on the initial data:
\begin{align}\label{iniz:vu}
   & \vu_0 \in V, \quad\dive \vu_0 = 0,\\
 \label{iniz:fhi}   
   & \fhi_0 \in H^3_{\per}(\Omega),\\
 \label{iniz:vt}
   & \vt_0 \in V, \quad \vt_0 \ge \ubarvt > 0~~\text{a.e.~in }\,\Omega,
\end{align}
where $\ubarvt$ is some positive constant (actually the last
condition in \eqref{iniz:vt} could be relaxed by 
asking $\vt_0>0$ almost everywhere
with $\log\vt\in L^1(\Omega)$).

\smallskip

With the above machinery at disposal, we can conclude this section
by stating our main existence theorem, whose proof will occupy the 
remainder of the paper:
\bete\label{teo:main}
 Let us assume\/ \eqref{hp:F1}-\eqref{hp:F3}, \eqref{hp:kappa},
 and\/ \eqref{iniz:vu}-\eqref{iniz:vt}. Let also $T>0$.
 Then, there exists at least
 one\/ {\rm strong solution} to the non-isothermal model for two-phase
 fluid flows, namely, one quadruple $(\vu,\fhi,\mu,\vt)$ with
 \begin{align} \label{rego:vu}
    & \vu \in H^1(0,T;H) \cap L^\infty(0,T;V) \cap L^2(0,T;H^2(\Omega)),\\
  \label{rego:fhi}
    & \fhi \in W^{1,\infty}(0,T;V') \cap H^1(0,T;V) \cap L^2(0,T;H^3(\Omega)),\\
   \label{rego:mu}
    & \mu \in H^1(0,T;V') \cap L^\infty(0,T;V) \cap L^2(0,T;H^3(\Omega)),\\
   \label{rego:vt}
    & \vt \in H^1(0,T;V') \cap L^\infty(0,T;L^{q+2}(\Omega)) \cap L^2(0,T;V),
      \quad \vt>0\quad\hbox{a.e. in }\,(0,T)\times \Omega,\\
   \label{rego:Kvt}
    & K(\vt) \in L^2(0,T;V),
 \end{align}  
 such that the equations of the system\/ \eqref{incom}-\eqref{CH2}
 hold in the sense of distributions as well as almost everywhere
 in~$(0,T)\times\Omega$, while~\eqref{calore} holds, for 
 a.e.~$t\in(0,T)$, as the following relation in~$V'$:
 \begin{equation}\label{calore2}
   \vt_t + \vu \cdot \nabla \vt 
    + \vt \big( \fhi_t + \vu \cdot \nabla \fhi \big) 
    - \Delta K(\vt) 
   = | \nabla \vu |^2 + | \nabla \mu |^2,
 \end{equation}
 where $\Delta$ is a weak form of the Laplace operator
 with periodic boundary conditions.
 Moreover, the quadruple $(\vu,\fhi,\mu,\vt)$ 
 complies with the initial condition
 \begin{equation}\label{iniz}
   \vu|_{t=0} = \vu_0, \qquad
   \fhi|_{t=0} = \fhi_0, \qquad
   \vt|_{t=0} = \vt_0,  
 \end{equation}
 almost everywhere in~$\Omega$.
\ente

%%%%%%%%%%%%%%%%%%%%%%%%%%%%%%%%%%%%%%%%%%%%%%%%%%%%%%%%%%%%%%%%%%%%%%%%%%%%%%%%%%%%%%%%%%%%%%%%%
%%%%%%%%%%%%%%%%%%%%%%%%%%%%%%%%%%%%%%%%%%%%%%%%%%%%%%%%%%%%%%%%%%%%%%%%%%%%%%%%%%%%%%%%%%%%%%%%%

\section{Global existence}
\label{sec:well}

We start by deriving the a-priori estimates leading to existence
of weak solutions. As noted in the Introduction, we shall work directly,
though formally, on the original system \eqref{incom}-\eqref{calore}
without referring to any approximation or regularization. Indeed, this
permits us to make the argument more readable and to avoid technical
complications. The details of a possible regularization scheme are posponed
to Section~\ref{sec:appro} below. In the following, the letter $c$ will
denote a generic positive constant depending only on the data of the
problem, whose value is allowed to vary on occurrence. In particular,
$c$ is intended to be independent of all regularization parameters.

%In the sequel (see \eqref{co75} or \eqref{co77d} for instance), we are noting by $p-$ (respectively, $p+$)
%an exponent strictly smaller (respectively, larger) than $p$, but which may be chosen arbitrarily close
%to $p$. Of course, the value of the constants $c$ is also influenced
%by the choice of~$p$. In particular $\infty-$ stands for an arbitrarily
%large positive number $P$ which can be chosen freely, in order to exploit
%the relation $\|v\|_{L^{\infty-}(\Omega)} \le c \|v\|_V$ coming from classical two-dimensional embeddings.

\smallskip

\noindent%
{\bf Energy estimate.}~~%
This basic property corresponds, in the physical derivation of the
model, to the energy conservation principle. To deduce it from the equations,
we test \eqref{ns} by $\vu$, \eqref{CH1} by $\mu$, \eqref{CH2} by $-\fhi_t$,
\eqref{calore} by $1$, integrate over $\Omega$, and sum all the obtained
relations together. Then, using the fact that 
\begin{equation}\label{co11} 
  (\vu \cdot \nabla \fhi) \mu 
    = (\vu \cdot \nabla \fhi) (- \Delta \fhi + F'(\fhi) - \vt )
\end{equation}
thanks to \eqref{CH2}, and performing standard integration by
parts (cf.~\cite[Sec.~2]{ERS}
for more details), it is not difficult to arrive at the relation
\begin{equation}\label{st:energy}
 \ddt \calE(\vu, \fhi, \vt) = 0, 
  \qquext{where }\,\calE(\vu, \fhi, \vt)
     := \io \bigg( \frac12 | \vu |^2 
        + \frac12 | \nabla \fhi |^2 
	+ F(\fhi)    
        + \vt \bigg)
\end{equation}
is the {\sl total energy}\/ of the system, given by the sum
of the {\sl kinetic}, {\sl interfacial}, {\sl configuration}, 
and {\sl thermal}\/ energies (the four summands in $\calE$).

Relation~\eqref{st:energy}, in turn, yields the following
{\sl a priori}\/ estimates: 
\begin{align}\label{st11}
  & \| \vu \|_{\LIH} \le c,\\
 \label{st12}
  & \| \fhi \|_{\LIV} \le c,\\ 
 \label{st13}
  & \| \vt \|_{L^\infty(0,T;L^1(\Omega))} \le c,
\end{align}  
where the control of the full $V$-norm of $\fhi$ (and not only
of the $L^2$-norm of the gradient) is reached thanks to the
superlinear growth of $F$ at infinity (cf.~\eqref{hp:F1}).
Note that, for getting \eqref{st13} from \eqref{st:energy},
the nonnegativity of $\vt$ is exploited
(which holds as a consequence of the 
approximation scheme, cf.~Lemma~\ref{lemma:fp3} below).
We also observe that, thanks to \eqref{st12} and Sobolev's 
embeddings, there follows
\begin{equation}\label{st14} 
  \| \fhi \|_{L^\infty(0,T;L^p(\Omega))} \le c_p
   \quext{for all }\,p\in[1,\infty).
\end{equation}

\smallskip

\noindent%
{\bf Conservation properties.}~~%
Integrating \eqref{ns} and \eqref{CH1} over $\Omega$, and
using \eqref{incom} together with
the periodic boundary conditions, we obtain
\begin{equation}\label{cons}
  \ddt \io \vu = \ddt \io \fhi = 0
   \quext{a.e.~in }\,(0,T).
\end{equation}
In other words, the spatial mean values of the velocity
and of the phase variable are constant in time. This basically
corresponds to the physical principles of conservation of
momentum and of mass. Of course, \eqref{cons} can
be equivalently rewritten as
\begin{equation}\label{zeromean}
  (\vu_t(t))\OO = (\fhi_t(t))\OO = 0
   \quext{for a.e.~}\,t\in(0,T).
\end{equation}

\smallskip

\noindent%
{\bf Entropy estimate.}~~%
The following estimate corresponds to the entropy production
principle. It is simply obtained by testing \eqref{calore} by
$-\vt^{-1}$ and integrating over~$\Omega$. As before, in order
for the procedure to be rigorous, we need that 
$\vt$ is strictly positive in the approximation.
Recalling also \eqref{hp:kappa}, we
then obtain
\begin{equation}\label{co21}
  \ddt \io ( - \log \vt - \fhi )
   + \io \frac1{\vt} \big( | \nabla \vu |^2 + | \nabla \mu |^2 \big)
   + \io \big( | \nabla \log \vt |^2 + k_q | \nabla \vt^{q/2} |^2 \big)
   = 0,
\end{equation}
where $k_q>0$ only depends on the exponent~$q$ 
(cf.~\eqref{hp:kappa}).
Integrating in time and recalling~\eqref{st12}-\eqref{st13},
we get the a priori bounds
\begin{align}\label{st21}
  & \| \log \vt \|_{L^\infty(0,T;L^1(\Omega))} 
    + \| \log \vt \|_{\LDV} \le c,\\
 \label{st22}
  & \| \nabla \vt^{q/2} \|_{\LDH} \le c.
\end{align}  
In particular, \eqref{st21} entails that the strict positivity
(almost everywhere in~$(0,T)\times \Omega$) of $\vt$ is preserved 
also in the limit. Note also that,
from \eqref{st22}, \eqref{st13} and inequality 
\eqref{poinc}, there follows
\begin{equation}\label{st23}
  \| \vt^{q/2} \|_{\LDV} \le c.
\end{equation}
Now, \eqref{st23} entails in particular
\begin{equation}\label{st24-a}
  \| \vt \|_{\LDH} \le c.
\end{equation}
On the other hand, arguing in a similar way as in~\cite[Sec. 4.2]{ERS}, 
by the elementary inequality
\begin{equation}\label{elemin}
  1 \le c \left (\frac{1}{x^2} + x^{q-2}\right ) \quext{for all }\,x\ge 0,
\end{equation}
holding for $q \ge 2$, we obtain from \eqref{st21}-\eqref{st22}
the additional bound
\begin{equation}\label{st24-b}
  \| \nabla \vt\|_{\LDH} \le c.
\end{equation}
Putting together \eqref{st24-a} and \eqref{st24-b}, we deduce 
\begin{equation}\label{st24}
  \| \vt \|_{\LDV} \le c.
\end{equation}

\smallskip

\noindent%
{\bf Temperature estimate.}~~%
While the bounds obtained above basically correspond to
physical principles, in order to get weak stability of families
of solutions we need to derive more refined a-priori estimates
by properly managing the equations of the system.
To this aim, we first integrate \eqref{calore} over~$\Omega$ with the
purpose of obtaining some information from the quadratic terms 
on the \rhs. Using the periodic boundary conditions, we 
actually infer
\begin{equation}\label{co31}
  \io \big( | \nabla \vu |^2 + | \nabla \mu |^2 \big)
   = \ddt \io \vt 
    + \io \vt ( \fhi_t + \vu \cdot \nabla \fhi )
\end{equation}
and we aim at controlling the terms on the \rhs. Actually, the
first one, after integration in time, is estimated simply
recalling~\eqref{st13}. Using \eqref{CH1} and H\"older's and 
Young's inequalities, we control the second integral as
follows:
\begin{equation}\label{co32}
  \io \vt ( \fhi_t + \vu \cdot \nabla \fhi )
   = \io \vt \Delta \mu
   = - \io \nabla \vt \cdot \nabla \mu
   \le \frac12 \big( \| \nabla \mu \|^2 + \| \nabla \vt \|^2 \big).
\end{equation}
The first term on the \rhs\ is absorbed
by the corresponding one on the \lhs\ of \eqref{co31},
while the latter is estimated thanks to~\eqref{st24}. 
Hence, we get
\begin{align}\label{st31}
  & \| \vu \|_{L^2(0,T;V)} \le c,\\
 \label{st32}
  & \| \nabla \mu \|_{\LDH} \le c.
\end{align}  
Now, integrating \eqref{CH2} in space, using \eqref{st13}, \eqref{st14} 
and \eqref{hp:F3}, and taking the (essential) supremum with 
respect to time, we readily infer
\begin{equation}\label{st33}
  \| \mu\OO \|_{L^\infty(0,T)} \le c.
\end{equation}
This property, combined with \eqref{st32}, yields
\begin{equation}\label{st34}
  \| \mu \|_{\LDV} \le c.
\end{equation}

\smallskip

\noindent%
The estimates proved up to this point basically correspond (up to the
slightly different assumptions on coefficients and data) to the procedure
used in~\cite{ERS} to get existence of at least one solution to 
a suitable reformulation of the problem. In particular, we recall
that, in the approach of~\cite{ERS}, the ``heat'' equation~\eqref{calore}
was restated in the form of a {\sl total energy balance}\/ complemented
with an {\sl energy production inequality}. 
Indeed, in that formulation the troublesome quadratic
terms characterizing the \rhs\ of \eqref{calore} do no longer appear.

It is worth noting that, in spite of the better embeddings
we have at disposal, the above estimates do not seem sufficient to prove
anything better in our 2D case, because, still, the \rhs\ of \eqref{calore}
is controlled only in $L^1$ (thanks to~\eqref{st31}-\eqref{st32}).
In view of these considerations, it is natural to investigate whether
higher order a-priori estimates could hold in the 2D case. 
As we will see in a while, such 
estimates can, indeed, be obtained, but the argument is not at all 
straighforward. For clarity, the procedure will be split into several
steps.
\beos\label{calolat}
 To be more precise, a further difference in our assumptions stands in the 
 choice of a {\sl linear}\/ latent heat (namely, we have no coefficient in
 front of $\vt_t$ in~\eqref{calore}), while in\/ \cite{ERS} we required
 to have a power-like function multiplying $\vt_t$, which yielded
 some additional summability of $\vt$. Actually, even in~2D, it is not
 clear whether one could prove existence of weak solutions (i.e., solutions
 complying with the formulation given in~\cite{ERS}) in the case of a linear
 latent heat. Nevertheless, this assumption is allowed as one looks
 for {\sl strong}\/ solutions, as we are doing now. 
\eddos     

\smallskip

\noindent%
{\bf Second estimate for $\fhi$.}~~%
We test \eqref{CH2} by $\Delta^2 \fhi$ and integrate over $\Omega$. 
Recalling \eqref{hp:F3}, we get
\begin{align}\no
   \| \nabla \Delta \fhi \|^2
  & = \io F''(\fhi) \nabla \fhi \cdot \nabla \Delta \fhi
    - \io \nabla ( \vt + \mu ) \cdot \nabla \Delta \fhi\\
 \no
  & \le c \| \nabla \Delta \fhi \| 
     \big( \| \nabla \fhi \| + \| \fhi \|_{L^{2p_F}(\Omega)}^{p_F}  
        \| \nabla \fhi \|_{L^\infty(\Omega)} 
      + \| \nabla \vt \| + \| \nabla \mu \| \big)\\
 \no
  & \le c \| \nabla \Delta \fhi \| 
     \big( 1 + \| \fhi \|_{V}^{1/2} \| \fhi \|_{H^3(\Omega)}^{1/2} 
      + \| \nabla \vt \| + \| \nabla \mu \| \big)\\
 \no
  & \le c \| \nabla \Delta \fhi \| 
     \big( 1 + \| \fhi \|_V
      + \| \fhi \|_{V}^{1/2} \| \nabla \Delta \fhi \|^{1/2} 
      + \| \nabla \vt \| + \| \nabla \mu \| \big)\\
 \label{co41}  
  & \le \frac12 \| \nabla \Delta \fhi \|^2
    + c \big( 1 + \| \nabla \vt \|^2 + \| \nabla \mu \|^2 \big).
\end{align}
Note that \eqref{dis:Linfty} 
(applied to $v = \nabla \fhi$) and
estimates~\eqref{st12},~\eqref{st14} 
have been used. Integrating \eqref{co41} in time 
and using \eqref{st24} and \eqref{st32}, we then obtain
\begin{equation}\label{st41}
  \| \fhi \|_{L^2(0,T;H^3(\Omega))} \le c.
\end{equation}
This property has some notable consequences. Firstly, 
testing \eqref{CH1} by nonzero $v\in V$, recalling also \eqref{st11},
we can notice that
\begin{align}\no
  \duav{\fhi_t,v} & = - \io \nabla \mu \cdot \nabla v
   - ( \vu \cdot \nabla \fhi, v )
    \le \| \nabla \mu \| \| \nabla v \| 
          + \| \vu \| \| \nabla \fhi \|_{L^\infty(\Omega)} \| v \| \\  
 \label{st41xx}
   & \le c \big( \| \nabla \mu \| 
         + \| \fhi \|_{H^3(\Omega)} \big)
     \| v \|_V.
\end{align}
Hence, dividing by $\| v \|_V$, passing to the supremum with 
respect to $v\in V\setminus\{0\}$,
squaring, integrating in time, and using \eqref{st34} and \eqref{st41}, we obtain
\begin{equation}\label{st34b}
  \| \fhi_t \|_{\LDVp} \le c.
\end{equation}
Moreover, \eqref{st41} permits us to get a useful bound for the 
last term in~\eqref{ns}. Indeed, using twice \eqref{dis:L4}, 
we have
\begin{align}\no
  \big\| \dive ( \nabla \fhi \otimes \nabla \fhi ) \big\|
   & \le c \| D^2 \fhi \|_{L^4(\Omega)} \| \nabla \fhi \|_{L^4(\Omega)} 
   \le c \| \fhi \|_{H^3(\Omega)}^{1/2} \| \fhi \|_{H^2(\Omega)} \| \fhi \|_V^{1/2}\\
 \label{co42}
   & \le c \| \fhi \|_{H^3(\Omega)} \| \fhi \|_V
     \le c \| \fhi \|_{H^3(\Omega)},
\end{align}
where the last inequality follows from \eqref{st12}. Then, squaring and
integrating over $(0,T)$, thanks to \eqref{st41} we infer
\begin{equation}\label{st42}
  \big\| \dive ( \nabla \fhi \otimes \nabla \fhi ) \big\|_{\LDH}
   \le c.
\end{equation}

\smallskip

\noindent%
{\bf Second estimate for $\vu$.}~~%
Property \eqref{st42} allows us to apply standard
regularity results to the 2D~Navier-Stokes system 
\eqref{incom}-\eqref{ns}, (see for instance~\cite[Sec.~9.6]{rob}), 
basically corresponding to testing \eqref{ns} by $-\Delta \vu$. 
As a consequence, we infer
\begin{equation}\label{st51}
  \| \vu \|_{\HUH} 
   + \| \vu \|_{\LIV} 
   + \| \vu \|_{\LDHD} 
   \le c.
\end{equation}

\smallskip

\noindent%
{\bf Key estimate: $\fhi$.}~~%
We start now with the key regularity estimate, which is 
obtained by combining in a suitable way equations \eqref{CH1},
\eqref{CH2} and \eqref{calore}.
At first, we deal with the Cahn-Hilliard
system. Namely, we take \eqref{CH1}, differentiate it 
with respect to time, and test the result by 
$\calN \fhi_t$. Correspondingly, we differentiate
\eqref{CH2} in time and test by $- \fhi_t$. Summing
the obtained relations, noting that a couple of terms
cancel in view of
\begin{equation}\label{co61}
  ( \Delta \mu_t, \calN \fhi_t )
   = - \big( ( -\Delta ) ( \mu_t - (\mu_t)\OO ), (- \Delta)^{-1} \fhi_t \big)
   = - \big( \mu_t - (\mu_t)\OO, \fhi_t \big)
   = - ( \mu_t, \fhi_t ),
\end{equation}
where we have used also \eqref{zeromean}, we then get
\begin{align}\no
  & \frac12 \ddt \| \fhi_t \|_{V'}^2
    + \| \nabla \fhi_t \|^2
    + \io (F''(\fhi)+\lambda) | \fhi_t |^2\\
  \label{co62}   
  & \mbox{}~~~~~
   = \lambda \| \fhi_t \|^2
   - \duav{ \vu_t \cdot \nabla \fhi, \calN \fhi_t }
   - \duav{ \vu \cdot \nabla \fhi_t, \calN \fhi_t }
   + ( \vt_t , \fhi_t ).
\end{align}
Thanks to~\eqref{hp:F2}, the last term on the \lhs\
is nonnegative. On the 
other hand, we need to control the \rhs. To this aim, we
first notice that, by \eqref{inter-gen} and \eqref{zeromean}, 
we can estimate the first term as
\begin{equation}\label{co63}
  \lambda \| \fhi_t \|^2
   \le \frac18 \| \nabla \fhi_t \|^2
    + c \| \fhi_t \|_{V'}^2.
\end{equation}
Next, using \eqref{dis:Linfty} and standard embeddings,
we infer
\begin{align}\no
  - \duav{ \vu_t \cdot \nabla \fhi, \calN \fhi_t }
   & = (\!( \vu_t \cdot \nabla \fhi, \fhi_t )\!)_{V_0'}
   \le \| \vu_t \cdot \nabla \fhi \|_{V'} \| \fhi_t \|_{V'}
   \le c \| \vu_t \cdot \nabla \fhi \| \| \fhi_t \|_{V'} \\
 \label{co64}   
   & \le c \| \vu_t \| \| \nabla \fhi \|_{L^\infty(\Omega)} \| \fhi_t \|_{V'}
    \le c \| \vu_t \|^2
     + c \| \fhi \|_{H^3(\Omega)}^2 \| \fhi_t \|_{V'}^2
\end{align}
and
\begin{align}\no
  - \duav{ \vu \cdot \nabla \fhi_t, \calN \fhi_t }
   & = (\!( \vu \cdot \nabla \fhi_t, \fhi_t )\!)_{V_0'}
   \le \| \vu \cdot \nabla \fhi_t \|_{V'} \| \fhi_t \|_{V'}
   \le c \| \vu \cdot \nabla \fhi_t \| \| \fhi_t \|_{V'} \\
 \label{co65}
   & \le c \| \vu \|_{L^\infty(\Omega)} \| \nabla \fhi_t \| \| \fhi_t \|_{V'}
    \le \frac18 \| \nabla \fhi_t \|^2
     + c \| \vu \|_{H^2(\Omega)}^2 \| \fhi_t \|_{V'}^2.
\end{align}
In the above formulas we noted by 
$(\!( \cdot, \cdot )\!)_{V_0'}$ the scalar product of $V_0'$
and used the fact that $(-\Delta)$ corresponds to
the Riesz operator from $V_0$ to $V_0'$.

Hence, on account of \eqref{co63}-\eqref{co65}, relation \eqref{co62}
takes the form
\begin{equation}\label{co62b}   
  \frac12 \ddt \| \fhi_t \|_{V'}^2
   + \frac34 \| \nabla \fhi_t \|^2
   + \io (F''(\fhi)+\lambda) | \fhi_t |^2
  \le M_1(t) \| \fhi_t \|_{V'}^2
   + M_2(t)  
   + ( \vt_t , \fhi_t ),
\end{equation}
where the functions 
\begin{equation}\label{co66}
  M_1(t) = c \big( 1 + \| \fhi \|_{H^3(\Omega)}^2
           + \| \vu \|_{H^2(\Omega)}^2  \big), \qquad
  M_2(t) = c \big( 1 + \| \vu_t \|^2 \big)
\end{equation}
lie (or, more precisely, are uniformly bounded w.r.t.~all approximation
parameters) in $L^1(0,T)$ thanks to \eqref{st41} and \eqref{st51}.
It remains to control the last term on the \rhs. This, however, 
requires to work with the energy equation~\eqref{calore}, which is
our next task.
     
\smallskip

\noindent%
{\bf Key estimate: $\vt$.}~~%
We start testing \eqref{calore} by $\fhi_t$ in order 
to compute the last term in \eqref{co62b}. We get
\begin{align}\no
  & ( \vt_t , \fhi_t ) 
    + \io \vt \fhi_t^2 
   = \io \vt \, \vu \cdot \nabla \fhi_t
    - \io \vt \, \fhi_t \, \vu \cdot \nabla \fhi \\
 \label{co71}   
  & \mbox{}~~~~~
   - \io \kappa(\vt) \nabla \vt \cdot \nabla \fhi_t
   + \io \big( | \nabla \vu |^2 + | \nabla \mu |^2 \big) \fhi_t,
\end{align}
where we have used \eqref{incom}. Let us provide an estimate for the terms on the
\rhs. Firstly, we have
\begin{equation}\label{co73}
  \io \vt \, \vu \cdot \nabla \fhi_t
   \le \| \vt \|_{L^4(\Omega)}
    \| \vu \|_{L^4(\Omega)} \| \nabla \fhi_t \|
   \le \frac1{16} \| \nabla \fhi_t \|^2
    + c \| \vt \|_{V}^2,
\end{equation}
%
%for small $\epsi_1>0$ to be chosen later, and correspondingly
%large $c_1$. 
%Note that the value $c_\epsi$ also
where we have also used estimate~\eqref{st51}. 

Next, by \eqref{inter-gen2} and H\"older's and Young's
inequalities, thanks also to 
\eqref{st12} and \eqref{st51}, we have 
\begin{align}\no
  - \io \vt \, \fhi_t \, \vu \cdot \nabla \fhi 
   & \le \| \vt \|_{L^{\infty-}(\Omega)}
     \| \fhi_t \|_{L^{2+}(\Omega)}
     \| \vu \|_{L^{\infty-}(\Omega)}
     \| \nabla \fhi \|
    \le c \| \vt \|_{V}
     \| \fhi_t \|_{L^{2+}(\Omega)}\\
 \label{co75} 
   & \le c \| \vt \|_{V}^2
     + c \| \fhi_t \|_{L^{2+}(\Omega)}^2
    \le c \| \vt \|_{V}^2 
     + \frac1{16}\| \nabla \fhi_t \|^2
     + c \| \fhi_t \|_{V'}^2.
\end{align}
Here (and below) the notation $\infty-$ stands for an exponent which
is chosen as large 
(i.e.~as close as infinity) as we need (in view of the
fact that $V\subset L^p(\Omega)$ for all $p\in [1,\infty)$). 
Correspondingly, $2+$ turns out to be larger than, but close to $2$.
Next,
\begin{equation}\label{co76}
  - \io \kappa(\vt) \nabla \vt \cdot \nabla \fhi_t
   \le 4 \io \kappa^2(\vt) | \nabla \vt |^2
   + \frac1{16} \| \nabla \fhi_t \|^2.
\end{equation}
Finally, using interpolation (cf.~\eqref{dis:L4} 
and \eqref{inter-gen}),
Young's inequality, and~\eqref{st51},
\begin{align}\no
  \io | \nabla \vu |^2 \fhi_t 
   & \le c \| \fhi_t \|_{L^{4}(\Omega)} 
    \| \nabla \vu \| \| \nabla \vu \|_{L^4(\Omega)} 
   \le c \| \fhi_t \|^{1/2} \| \nabla \fhi_t \|^{1/2} 
    \| \vu \|_{H^2(\Omega)} \\
 \label{co79a}
   & \le c \| \fhi_t \|_{V'}^{1/4}
    \| \nabla \fhi_t \|^{3/4}
    \| \vu \|_{H^2(\Omega)}
   \le c + c \| \vu \|_{H^2(\Omega)}^2 
   + c \| \fhi_t \|_{V'}^2
   + \frac1{16} \| \nabla \fhi_t \|^2.
\end{align}
Thanks to \eqref{co73}-\eqref{co79a}, \eqref{co71}
gives
\begin{align}\no
  & ( \vt_t , \fhi_t ) 
    + \io \vt \fhi_t^2 
   \le c \big( 1 
    + \| \vu \|_{H^2(\Omega)}^2 
    + \| \vt \|_{V}^2
    + \| \fhi_t \|_{V'}^2 \big)\\
 \label{co71new}
  & \mbox{}~~~~~
   + \frac14 \| \nabla \fhi_t \|^2
   + 4 \io \kappa^2(\vt) | \nabla \vt |^2
   + \io | \nabla \mu |^2 \fhi_t.
\end{align}
The last two terms on the \rhs\ are still 
to be controlled. The quadratic term in $\nabla \mu$,
which is the most difficult one, will be dealt with
at the end. In order to treat the term in $\nabla \vt$,
we need another estimate. Namely, we test \eqref{calore}
by $8 K(\vt)$ (cf.~\eqref{defiK}). Let us set
\begin{equation}\label{defiH}
  \mathcal{J}(r) := \int_0^r K(s)\,\dis
     = \frac{r^2}2 + \frac{1}{(q+1)(q+2)} r^{q+2}, \quad
     r \ge 0.
\end{equation}
Then, we obtain
\begin{align}\no
  & 8 \ddt \io \mathcal{J}(\vt)
   + 8 \io \vu \cdot \nabla \mathcal{J}(\vt)
    + 8 \io \kappa^2(\vt) | \nabla \vt |^2\\
 \label{co72}
  & \mbox{}~~~~~
   = - 8 \io \vt K(\vt) ( \fhi_t + \vu \cdot \nabla \fhi )
   + 8 \io K(\vt) \big( | \nabla \vu |^2 + | \nabla \mu |^2 \big),
\end{align}
where, in fact, the second integral on the \lhs\ is zero
in view of \eqref{incom} and the periodic boundary 
conditions. Again, we need to provide an estimate for the
terms on the \rhs. At first, recalling \eqref{defiK}
and using \eqref{inter-gen}, we have
\begin{align}\no
  - 8 \io \vt K(\vt) \fhi_t
   & \le c \io \big( \vt^2 + \vt^{q+2} \big) | \fhi_t | 
     \le c \io \big( 1 + \vt^{q+2} \big) | \fhi_t | \\
 \label{co77}
   & \le c \| \fhi_t \| + \io \vt^{q+2} | \fhi_t | \le c + \frac1{16} \| \nabla \fhi_t \|^2
   + c \| \fhi_t \|_{V'}^2
   + c \io \vt^{q+2} | \fhi_t |.
\end{align}
The last term needs to be managed accurately. We start by
noting that
\begin{equation}\label{co77b}
  c \io \vt^{q+2} | \fhi_t |
   \le c \| \vt^{q+2} \| \| \fhi_t \|
   \le c \| \vt \|_{L^{2q+4}(\Omega)}^{q+2} 
       \| \fhi_t \|_{V'}^{1/2} \| \nabla \fhi_t \|^{1/2}.
\end{equation}
Next, we observe that, using \eqref{dis:Lp} with the choices $r = 2q+4$ 
and $s = 1$, we would obtain
\begin{equation}\label{co77c}
  \| \vt \|_{L^{2q+4}(\Omega)} 
   \le \| \vt \|_{L^1(\Omega)}^{\frac1{2q+4}}
         \| \vt \|_{L^\infty(\Omega)}^{\frac{2q+3}{2q+4}}.
\end{equation}
However the above interpolation exponents do not work
in our case and need to be modified slightly. Actually, recalling 
also \eqref{st13}, 
we can continuate \eqref{co77b} as follows:
\begin{align}\no
  c \io \vt^{q+2} | \fhi_t |
   & \le c \Big( \| \vt \|_{L^{1}(\Omega)}^{\frac1{2q+4}-}
    \| \vt \|_{L^{\infty-}(\Omega)}^{\frac{2q+3}{2q+4}+} \Big)^{q+2}
       \| \fhi_t \|_{V'}^{1/2} \| \nabla \fhi_t \|^{1/2} \\
 \label{co77d}  
   & \le c \| \vt \|_{L^{\infty-}(\Omega)}^{\frac{2q+3}{2}+}
       \| \fhi_t \|_{V'}^{1/2} \| \nabla \fhi_t \|^{1/2}  
   \le c \| \vt^{q+1} \|_{L^{\infty-}(\Omega)}^{\frac{2q+3}{2q+2}+}
       \| \fhi_t \|_{V'}^{1/2} \| \nabla \fhi_t \|^{1/2}.
\end{align}
As above, $\infty-$ stands for an exponent $P \in [1,\infty)$
(which we can choose as large as we need) and the number
$\frac{2q+3}{2}+$ depends on the choice of $P$
and will be closer to $\frac{2q+3}{2}$ as larger is taken $P$.
The same applies to $\frac{2q+3}{2q+2}+$ and other exponents below. 
Of course, also the constants $c$ 
will depend on the choice of $P$ (and will be larger for
larger $P$). Then, recalling \eqref{poinc}
and subsequently using Young's inequality with 
exponents~$8$, $4$, and $8/5$, computation \eqref{co77d}
can be continued this way:
\begin{align}\no
  c \io \vt^{q+2} | \fhi_t |
   & \le c \| \fhi_t \|_{V'}^{1/2} \| \nabla \fhi_t \|^{1/2}
    \big( 1 + \| \nabla \vt^{q+1} \|^{\frac{2q+3}{2q+2}+} \big)\\
 \label{co77e}  
   & \le c + c \| \fhi_t \|_{V'}^4
    + \frac1{16} \| \nabla \fhi_t \|^2 
    + c \| \nabla \vt^{q+1} \|^{\frac{4(2q+3)}{5(q+1)}+},
\end{align}
and, thanks to \eqref{hp:kappa} (actually, $q>1$ would be enough
at this level), we can take $P$ so large
that $\frac{4(2q+3)}{5(q+1)}+$ is {\sl strictly}\/ smaller
than $2$. Hence, using Young's inequality again, and recalling
\eqref{K11}, we conclude that
\begin{align}\no
  - 8 \io \vt K(\vt) \fhi_t
   & \le c + \frac18 \| \nabla \fhi_t \|^2 
    + c \| \fhi_t \|_{V'}^2 
    + c \| \fhi_t \|_{V'}^4 
    + c \| \nabla \vt^{q+1} \|^{\frac{4(2q+3)}{5(q+1)}+} \\
 \label{co77f}
  & \le c + \frac18 \| \nabla \fhi_t \|^2 
    + c \| \fhi_t \|_{V'}^4 
    + \io \kappa^2(\vt) | \nabla \vt |^2.
\end{align}
The estimation of the subsequent summand in \eqref{co72} is 
simpler. Actually,
recalling also \eqref{st12}, \eqref{st13} and \eqref{st51},
and using once more \eqref{poinc}, 
\eqref{K11} and Young's inequality, we get
\begin{align}\no
  - 8 \io \vt K(\vt) \vu \cdot \nabla \fhi
   & \le c \io ( 1 + \vt^{q+2} ) | \vu | \, | \nabla \fhi |
   \le c \big( 1 + \| \vt^{q+2} \|_{L^{2+}(\Omega)} \big)
     \| \vu \|_{L^{\infty-}(\Omega)} \| \nabla \fhi \| \\
 \label{co78a}
   & \le c \Big( 1 
        + \| \vt^{q+1} \|_{L^{\frac{2(q+2)}{q+1}+}(\Omega)}^{\frac{q+2}{q+1}} \Big)
    \le c \big ( 1 + \| \nabla \vt^{q+1} \|^{\frac{q+2}{q+1}} \big)
    \le c + \io \kappa^2(\vt) | \nabla \vt |^2,
\end{align}
since, clearly, $\frac{q+2}{q+1} < 2$. Next, 
recalling \eqref{st51} and \eqref{stima:Kvt} and arguing 
analogously to \eqref{co79a}, we have
\begin{equation}\label{co79b}
  8 \io | \nabla \vu |^2 K(\vt)
   \le c \| K(\vt) \|_{L^{4}(\Omega)} 
    \| \nabla \vu \| \| \nabla \vu \|_{L^4(\Omega)} 
   \le c + c \| \vu \|_{H^2(\Omega)}^2 
   + \io \kappa^2(\vt) | \nabla \vt |^2.
\end{equation}
Collecting \eqref{co77}-\eqref{co79b}, \eqref{co72}
gives
\begin{align}\no
  & 8 \ddt \io \mathcal{J}(\vt)
   + 5 \io \kappa^2(\vt) | \nabla \vt |^2\\
 \label{co72new}
  & \mbox{}~~~~~
   \le c + c \| \vu \|_{H^2(\Omega)}^2 
    + \frac18 \| \nabla \fhi_t \|^2 
    + c \| \fhi_t \|_{V'}^4 
    + 8 \io K(\vt) | \nabla \mu |^2.
\end{align}
Summing \eqref{co71new} and \eqref{co72new} we then get
\begin{align}\no
  & ( \vt_t , \fhi_t ) 
    + \io \vt \fhi_t^2 
    + 8 \ddt \io \mathcal{J}(\vt)
    + \io \kappa^2(\vt) | \nabla \vt |^2\\
 \label{co71+72}
  & \mbox{}~~~~~
   \le c \big( 1 
    + \| \vu \|_{H^2(\Omega)}^2 
    + \| \vt \|_{V}^2
    + \| \fhi_t \|_{V'}^4 \big)
    + \frac38 \| \nabla \fhi_t \|^2
   + \io ( 8 K(\vt) + \fhi_t ) | \nabla \mu |^2.
\end{align}
Hence, summing \eqref{co62b} and \eqref{co71+72} we obtain
\begin{align}\no
  & \frac12 \ddt \| \fhi_t \|_{V'}^2
   + 8 \ddt \io \mathcal{J}(\vt)
   + \frac38 \| \nabla \fhi_t \|^2
   + \io (F''(\fhi)+\lambda) | \fhi_t |^2\\
 \no  
  & \mbox{}~~~~~ 
   + \io \vt \fhi_t^2 
   + \io \kappa^2(\vt) | \nabla \vt |^2
  \le M_1(t) \| \fhi_t \|_{V'}^2
   + M_2(t) \\
 \label{key1old}  
  & \mbox{}~~~~~ 
   + c \big( 1 + \| \vu \|_{H^2(\Omega)}^2 
    + \| \vt \|_{V}^2
    + \| \fhi_t \|_{V'}^4 \big)
    + \io ( 8 K(\vt) + \fhi_t ) | \nabla \mu |^2.
\end{align}
Neglecting some positive terms in the \lhs\ and
rearranging, we then arrive at
\begin{align}\no
  & \frac12 \ddt \| \fhi_t \|_{V'}^2
   + 8 \ddt \io \mathcal{J}(\vt)
   + \frac38 \| \nabla \fhi_t \|^2
   + \io \kappa^2(\vt) | \nabla \vt |^2\\
  \label{key1}  
  & \mbox{}~~~~~ 
    \le M_3(t) \| \fhi_t \|_{V'}^2
   + M_4(t) 
   + \io ( 8 K(\vt) + \fhi_t ) | \nabla \mu |^2,
\end{align}
where we have set
\begin{equation}\label{co66new}
  M_3(t) = c \big( 1 + \| \fhi \|_{H^3(\Omega)}^2
           + \| \vu \|_{H^2(\Omega)}^2 
	   + \| \fhi_t \|_{V'}^2 \big), \qquad
  M_4(t) = c \big( 1 + \| \vu_t \|^2 + \| \vu \|_{H^2(\Omega)}^2 
    + \| \vt \|_{V}^2\big).
\end{equation}

\smallskip

\noindent%
{\bf Key estimate: quadratic terms.}~~%
The most difficult part of our argument concerns the control of 
the last term in the \rhs\ of \eqref{key1}. This is 
based on the embedding inequality 
\eqref{yudo} proved in Lemma~\ref{lemma:yudo} below. 
Indeed, applying \eqref{yudo} to $\xi = | \nabla \mu |^2$
and using once more \eqref{stima:Kvt} together with \eqref{zeromean} 
and the Poincar\'e-Wirtinger inequality, we get
\begin{align}\no
  & \io ( 8 K(\vt) + \fhi_t ) | \nabla \mu |^2
   \le c \big( \| K(\vt) \|_V + \| \nabla \fhi_t \| \big)
         \big \| | \nabla \mu |^2 \|_{V'}\\
 \no	 
   & \mbox{}~~~~~ \le c + \frac12 \io \kappa^2(\vt) | \nabla \vt |^2
    + \frac18 \| \nabla \fhi_t \|^2
     + c \big \| | \nabla \mu |^2 \|_{V'}^2\\
 \label{key31}
  & \mbox{}~~~~~ \le c + \frac12 \io \kappa^2(\vt) | \nabla \vt |^2
    + \frac18 \| \nabla \fhi_t \|^2
    + c \big\| | \nabla \mu |^2 \big\|_{L^1(\Omega)}^2 
          \log \big(e + \big\| | \nabla \mu |^2 \big\|_{L^2(\Omega)} \big).
\end{align}
Now, let us notice that $\psi(r) = e^r$, $r\in \RR$ and
$\psi^*(s) = s (\log s - 1)$, $s > 0$ (extended by continuity
to $s=0$ by setting $\psi^*(0)=0$) are {\sl convex conjugate}\/
functions. Consequently, for any $r\in \RR$, $s\ge 0$, we 
have (cf., e.g., \cite[Sec.~1.4]{BrAF})
$rs \le \psi(r) + \psi^*(s)$.
Applying this property 
to~$r = \log \big( e + \| | \nabla \mu |^2 \|_{L^2(\Omega)} \big)$
and $s = c \big\| | \nabla \mu |^2 \big\|_{L^1(\Omega)}^2$,
the last term (let us note it as $\calI$) in \eqref{key31} 
can be controlled as follows:
\begin{align}\no
  \calI & \le c \big\| | \nabla \mu |^2 \big\|_{L^1(\Omega)}^2 
   \Big( \log \big( c \big\| | \nabla \mu |^2 \big\|_{L^1(\Omega)}^2 \big) - 1 \Big)
    + e + \big\| | \nabla \mu |^2 \big\|_{L^2(\Omega)}\\
 \label{key32}
   & \le c + c \| \nabla \mu \|^4
   \log \big( e + \| \nabla \mu \|^2 \big)
     + \| \nabla \mu \|_{L^4(\Omega)}^2,
\end{align}
where, observing that 
$\big\| | \nabla \mu |^2 \big\|_{L^1(\Omega)}^2  = \|\nabla \mu\|^4$, 
we used the fact that
\begin{align}\no
 \big\| | \nabla \mu |^2 \big\|_{L^1(\Omega)}^2  \log \big (c \big\| | \nabla \mu |^2 \big\|_{L^1(\Omega)}^2  \big ) 
  & = \|\nabla \mu\|^4 \log \big (c \|\nabla \mu\|^4) \\
 \no & = \|\nabla \mu\|^4 \big (\log c + 2 \log \|\nabla \mu\|^2 \big )\\
 \no 
  & \le c \|\nabla \mu\|^4 \big (1 + \log \|\nabla \mu\|^2 \big ) \\
 \label{key32bis}
 & \le c \|\nabla \mu\|^4 \log \big(e + \|\nabla \mu\|^2\big );
\end{align}
the last line is based on the elementary inequality 
$1 + \log \lambda \le c \log (e + \lambda)$, holding for $\lambda > 0$.
Now, to manage the last term of \eqref{key32}, we use equation~\eqref{CH1},
estimate~\eqref{st12}, and inequalities~\eqref{dis:L4} and \eqref{inter-gen}:
\begin{align}\no
  \| \nabla \mu \|_{L^4(\Omega)}^2 
   & \le c \| \nabla \mu \| \| \mu \|_{H^2(\Omega)}
    \le c \| \nabla \mu \| \big( \| \mu \|_{V} + \| \Delta \mu \| \big)
    \le c \| \mu \|_{V}^2 + c \| \fhi_t \|^2 + c \| \vu \cdot \nabla \fhi \|^2 \\
 \no
   & \le c \| \mu \|_{V}^2 + \frac18 \| \nabla \fhi_t \|^2 
    + c \| \fhi_t \|_{V'}^2 
    + c \| \vu \|_{L^\infty(\Omega)}^2 \| \nabla \fhi \|^2 \\
 \label{key33}
  & \le c \| \mu \|_{V}^2 + \frac18 \| \nabla \fhi_t \|^2 
    + c \| \fhi_t \|_{V'}^2 + c \| \vu \|_{H^2(\Omega)}^2. 
\end{align}
On the other hand, the first nonconstant term on the \rhs\ 
of \eqref{key32} needs a further manipulation. 
Namely, we have to test equation \eqref{CH1} by 
$- \mu$. Then, noting that both terms on the \lhs\ of \eqref{CH1}
have zero spatial mean and using the Poincar\'e-Wirtinger
inequality, we get
\begin{align}\no
  \| \nabla \mu \|^2 
   & = - \io ( \fhi_t + \vu \cdot \nabla \fhi ) \mu 
    = - \io ( \fhi_t + \vu \cdot \nabla \fhi ) ( \mu - \mu\OO )
    \le \| \mu - \mu\OO \|_V \| \fhi_t + \vu \cdot \nabla \fhi \|_{V'}\\
 \no
   & \le \frac12 \| \nabla \mu \|^2
    + c \| \fhi_t \|^2_{V'} + c \| \vu \cdot \nabla \fhi \|^2_{V'} 
    \le \frac12 \| \nabla \mu \|^2
    + c \| \fhi_t \|_{V'}^2 + c \| \vu \cdot \nabla \fhi \|_{L^{4/3}(\Omega)}^2 \\
 \label{key34}
   & \le \frac12 \| \nabla \mu \|^2
    + c \| \fhi_t \|_{V'}^2 + c \| \vu \|_{L^4(\Omega)}^2 \| \nabla \fhi \|^2 
   \le \frac12 \| \nabla \mu \|^2
    + c \| \fhi_t \|_{V'}^2 + c,
\end{align}
where in the last line we also used \eqref{st51}. Consequently,
\begin{align}\no
  c \| \nabla \mu \|^4 \log \big( e + \| \nabla \mu \|^2 \big)
   & \le c \big( 1 + \| \fhi_t \|_{V'}^4 \big)
    \log \Big( e + c \big( 1 + \| \fhi_t \|_{V'}^2 \big) \Big)\\
 \label{key35}   
   & \le c + c \| \fhi_t \|_{V'}^4 
    \log \big( e + \| \fhi_t \|_{V'}^2 \big).
\end{align}
Hence, collecting \eqref{key32}-\eqref{key35}, 
\eqref{key31} gives
\begin{align}\no
  & \io ( 8 K(\vt) + \fhi_t ) | \nabla \mu |^2
   \le c + \frac12 \io \kappa^2(\vt) | \nabla \vt |^2
    + \frac14 \| \nabla \fhi_t \|^2
    + c \| \fhi_t \|^2_{V'}
    + c \| \mu \|^2_{V}\\
 \label{key36}   
 & \mbox{}~~~~~   
   + c \| \vu \|_{H^2(\Omega)}^2 
   + c \| \fhi_t \|_{V'}^4 
    \log \big( e + \| \fhi_t \|_{V'}^2 \big).
\end{align}
Plugging \eqref{key36} into \eqref{key1} we then get
\begin{align}\no
  & \frac12 \ddt \| \fhi_t \|_{V'}^2
   + 8 \ddt \io \mathcal{J}(\vt)
   + \frac18 \| \nabla \fhi_t \|^2
   + \frac12 \io \kappa^2(\vt) | \nabla \vt |^2\\
  \label{key2}  
  & \mbox{}~~~~~ 
    \le M_5(t) \| \fhi_t \|_{V'}^2
   + M_6(t) 
   + c \| \fhi_t \|_{V'}^4 
    \log \big( e + \| \fhi_t \|_{V'}^2 \big).
\end{align}
where
\begin{align}\label{co66z1}
  M_5(t) & = c \big( 1 + \| \fhi \|_{H^3(\Omega)}^2
           + \| \vu \|_{H^2(\Omega)}^2 
	   + \| \fhi_t \|_{V'}^2 \big), \\
 \label{co66z2}
  M_6(t) & = c \big( 1 + \| \vu_t \|^2 + \| \vu \|_{H^2(\Omega)}^2 
    + \| \vt \|_{V}^2 + \| \mu \|_V^2 \big).
\end{align}
Let us now set
\begin{equation}\label{defiY}
  Y_1(t):= \frac12 \| \fhi_t(t) \|_{V'}^2, \qquad
   Y_2(t):= 8 \io \mathcal{J}(\vt(t)).
\end{equation}
Then, from \eqref{key36} we obtain the following
differential inequality:
\begin{equation}\label{ineqY}
   Y_1'(t) + Y_2'(t)  \le M_5(t) Y_1(t)
   + M_6(t) + c Y_1^2(t)
    \log \big( e + Y_1(t) \big).
\end{equation}
Setting $Z(t) := e + Y_1(t) + Y_2(t)$, and dividing both hand sides
of \eqref{ineqY} by $Z \log Z$, it is then easy to get
\begin{equation}\label{ineqZ}
   \ddt \log \log Z(t)
   = \frac{Z'(t)}{Z(t) \log Z(t)}  
   \le M_5(t) + \frac{M_6(t)}{Z(t) \log Z(t)}
    + c Y_1(t).
\end{equation}
Now, let us notice that, in view of the a-priori estimates 
\eqref{st23}, \eqref{st24}, \eqref{st34}, \eqref{st41},
\eqref{st34b}, and \eqref{st51}, we have
\begin{equation}\label{YL1}
  \| Y_1 \|_{L^1(0,T)}
   + \| M_5 \|_{L^1(0,T)} 
   + \| M_6 \|_{L^1(0,T)}
   \le c.
\end{equation}
Moreover, we can notice that, at least formally,
\begin{align}\no
  Z(0) & = e + \frac12 \| \fhi_t(0) \|_{V'}^2
   + 8 \io \mathcal{J}(\vt_0) \\
 \no  
   & \le c + c \| \Delta \mu(0) \|_{V'}^2 + c \| \vu_0 \cdot \nabla \fhi_0 \|_{V'}^2
   + c \io \mathcal{J}(\vt_0) \\
 \no  
   & \le c + \frac12 \| \mu(0) \|_{V}^2 + c \| \vu_0 \|_V^2 \| \fhi_0 \|_{H^3(\Omega)}^2
   + c \| \vt_0 \|^2 
   + c \| \vt_0 \|_{L^{q+2}(\Omega)}^{q+2}\\
 \label{Z00}
   & \le c + c \| \Delta\fhi_0 \|_{V}^2 
    + c \| F'(\fhi_0) \|_{V}^2
    + c \| \vt_0 \|_{V}^{q+2}
    < \infty,
\end{align}
where in the second row we used equation~\eqref{CH1},
in the third we used standard interpolation and embeddings
and the definition of $\mathcal{J}$, and in the fourth we used 
equation \eqref{CH2} and the assumptions \eqref{hp:F3} on
the potential and \eqref{iniz} on the initial data.

Hence, thanks to \eqref{YL1} and to \eqref{Z00},
we can integrate \eqref{ineqZ} over $(0,T)$
to obtain
\begin{equation}\label{ZLinfty}
  \| Z \|_{L^\infty(0,T)} \le c.
\end{equation}
In particular, as often happens in 2D models, $Z$ turns out to 
grow with respect to time (at most)
as fast as a double exponential; however it
does not explode in finite times. This is the key point
in our existence proof.

\smallskip

\noindent%
{\bf Key estimate: consequences.}~~%
Recalling also \eqref{defiH}, \eqref{ZLinfty} gives
\begin{align}\label{stkey1}
  & \| \fhi_t \|_{\LIVp} \le c,\\
 \label{stkey2}
  & \| \vt \|_{L^\infty(0,T;L^{q+2}(\Omega))} \le c.
\end{align}  
Using \eqref{CH1} we can estimate, as in \eqref{key34}, the 
$H$-norm of $\nabla \mu$ in terms of the $V'$-norm
of $\fhi_t$. Hence, recalling also \eqref{st33}, we get
\begin{equation}\label{stkey3}
  \| \mu \|_{\LIV} \le c.
\end{equation}
With these properties at disposal, we can go back 
to relation~\eqref{key2}. Integrating it over $(0,T)$, 
we infer (cf.~also \eqref{defiK}, \eqref{zeromean})
\begin{align}\label{stkey4}
  & \| \fhi_t \|_{\LDV} \le c,\\
 \label{stkey5}
  & \| K(\vt) \|_{\LDV} \le c.
\end{align}  
In addition to that, an easy interpolation argument
and estimates \eqref{st41}, \eqref{st51}
permit us to check that
\begin{equation}\label{co88}
  \| \vu \cdot \nabla \fhi \|_{\LDV}
   \le c \| \vu \|_{\LIV} \| \fhi \|_{L^2(0,T;H^3(\Omega))}
   \le c.
\end{equation}
Hence, viewing \eqref{CH1} as a time-dependent family of elliptic
problems and using standard regularity results with \eqref{co88} and
\eqref{stkey4}, we infer
\begin{equation}\label{stkey6}
  \| \mu \|_{L^2(0,T;H^3(\Omega))} \le c.
\end{equation}
\beos\label{inizprimo}
 The argument given in \eqref{Z00} in order
 to estimate $\fhi_t(0)$ in terms
 of $\vu_0$, $\vt_0$ and $\mu_0$ is just formal since,
 at least in principle, it is not clear whether 
 equation~\eqref{CH1} holds pointwise for {\sl all}\/ values
 of the time variable. However, as will be clear from 
 the approximation scheme detailed in the next section,
 at least locally in time we can get an approximate solution
 which is essentially as smooth as we need. Hence, at the 
 price of some additional technicalities, the argument could be 
 made fully rigorous with a limited effort.
\eddos

\smallskip
 
\noindent%
{\bf Estimate of $\mu_t$ and $\vt_t$.}~~%
In order to pass to the limit in the nonlinear terms involving
$\vt$ and $\mu$ we will apply the Aubin-Lions lemma. To this
aim, we need to deduce some a-priori estimates on $\mu_t$
and $\vt_t$. Hence, we (formally) differentiate
\eqref{CH2} with respect to time and use \eqref{calore} 
to get
\begin{equation}\label{CH2calore}
  \mu_t = - \Delta \fhi_t + F''(\fhi)\fhi_t
    + \vu \cdot \nabla \vt + \vt (\fhi_t + \vu \cdot \nabla \fhi)
    - \Delta K(\vt) - | \nabla \vu |^2 - | \nabla \mu |^2.
\end{equation}
Here, to justify the procedure, one could make similar observations as
in Remark~\ref{inizprimo} (see also Remark~\ref{remfp3} below).
Let us test the above relation by nonzero $v\in V$. We get,
using our choice of boundary conditions,
\begin{align}\no
  \duav{ \mu_t , v} & = \io \nabla ( \fhi_t + K(\vt) )\cdot \nabla v
   + \big( F''(\fhi) \fhi_t + \vu \cdot \nabla \vt 
    + \vt (\fhi_t + \vu \cdot \nabla \fhi) - | \nabla \vu |^2 - | \nabla \mu |^2 , v \big)\\
 \no
  & \le \| v \|_V \Big( \big\| \nabla ( \fhi_t + K(\vt) ) \big\| \\
 \label{coa1}   
  & \mbox{}~~~~~~~~~~
   + \big\| F''(\fhi) \fhi_t + \vu \cdot \nabla \vt 
    + \vt (\fhi_t + \vu \cdot \nabla \fhi) - | \nabla \vu |^2  
        - | \nabla \mu |^2 \big\|_{L^{3/2}(\Omega)} \Big).
\end{align}
Dividing by $\| v \|$, passing to the supremum with respect to
$v\in V\setminus\{0\}$, squaring, and integrating in time,
we would then obtain
\begin{equation}\label{stmut}
  \| \mu \|_{\HUVp} \le c,
\end{equation}
provided we could prove that
\begin{align}\label{cond11}
  & \big\| \nabla ( \fhi_t + K(\vt) ) \big\|_{\LDH} \le c,\\
 \label{cond12} 
  & \big \| F''(\fhi) \fhi_t + \vu \cdot \nabla \vt 
    + \vt (\fhi_t + \vu \cdot \nabla \fhi) - | \nabla \vu |^2 
   - | \nabla \mu |^2 \big\|_{L^2(0,T;L^{3/2}(\Omega))} \le c,
\end{align}
where the exponent $3/2$ is chosen just for simplicity
(any number strictly greater than~$1$ would be allowed, indeed).
Now, \eqref{cond11} is an immediate consequence of \eqref{stkey4}-\eqref{stkey5},
whereas \eqref{cond12} follows by appropriately combining all the above
a-priori estimates and using standard inequalities.
The details are lengthy but straighforward; hence, they 
are left to the reader.

To conclude, we test equation \eqref{calore} by $v\in V\setminus\{0\}$. 
Performing the very same computations as above we then get
\begin{equation}\label{stvtt}
  \| \vt \|_{\HUVp} \le c.
\end{equation}
This is the last estimate we need.

\medskip
 
\noindent%
{\bf Weak sequential stability.}~~%
We assume to have a sequence of solutions $(\vu_n,\fhi_n,\mu_n,\vt_n)$
satisfying the a-priori estimates obtained above uniformly with
respect to $n$. This could be, for instance, a sequence
of approximate solutions provided by the fixed-point argument
performed in the next section. Then, we aim at proving
that, up to the extraction of a subsequence, we can find
a limit quadruple $(\vu,\fhi,\mu,\vt)$ satisfying
\eqref{incom}-\eqref{calore} in the sense of 
Theorem~\ref{teo:main}. Note that, in view of the uniform 
character of the estimates, even though the approximate
solutions are defined only locally in time, we will have 
global solutions in the limit (cf.~also Remark~\ref{remfp1}
below). For this reason, and also for the sake of 
simplicity, we shall directly work on the 
original reference time interval~$(0,T)$.
 
That said, the above a-priori estimates  
(cf., in particular, 
\eqref{st24}, \eqref{st41}, \eqref{st51},
\eqref{stkey1}-\eqref{stkey6}, \eqref{stmut},
and \eqref{stvtt}), together with standard 
weak compactness results, entail 
\begin{align}\label{covu}
  & \vu_n \to \vu
   \quext{weakly star in }\, \HUH 
     \cap L^\infty(0,T;V) \cap L^2(0,T;H^2(\Omega)),\\
 \label{cofhi}
  & \fhi_n \to \fhi 
   \quext{weakly star in }\, W^{1,\infty}(0,T;V') 
     \cap H^1(0,T;V) \cap L^2(0,T;H^3(\Omega)),\\
 \label{comu}
  & \mu_n \to \mu
   \quext{weakly star in }\, \HUVp \cap L^\infty(0,T;V) 
     \cap L^2(0,T;H^3(\Omega)),\\
 \label{covt}
  & \vt_n \to \vt
   \quext{weakly star in }\, \HUVp \cap L^\infty(0,T;L^{q+2}(\Omega))
     \cap L^2(0,T;V).
\end{align}
Here and below, all convergence relations are intended to
hold up to the extraction of non-relabelled subsequences.
Note also that \eqref{cofhi}, by interpolation 
(cf., e.g., \cite[formula~(2.5.38)]{BG}), yields
\begin{equation}\label{cofhi2}
  \fhi_n \to \fhi 
   \quext{weakly in }\, H^{s}(0,T;H^{3-2s}(\Omega))
   \quext{for all }\,s\in [0,1].
\end{equation}
Then, the Aubin-Lions compactness theorem 
(cf.~\cite[p.~58]{lions} or~\cite[Cor.~4, Sec.~8]{simon}) 
entails (actually something more is true)
\begin{align}\label{covus}
  & \vu_n \to \vu
   \quext{strongly in }\, C^0([0,T];H) 
    \cap L^2(0,T;V),\\
 \label{cofhis}
  & \fhi_n \to \fhi 
   \quext{strongly in }\, C^0([0,T];V) 
    \cap L^2(0,T;H^2(\Omega)),\\
  \label{comus}
  & \mu_n \to \mu
   \quext{strongly in }\, C^0([0,T];H) 
    \cap L^2(0,T;H^2(\Omega)),\\
 \label{covts}
  & \vt_n \to \vt \quext{strongly in }\, C^0([0,T];V') 
    \cap \LDH.
\end{align}
Notice that, to deduce the first \eqref{cofhis}, 
also \eqref{cofhi2} has been used together with the
compact embedding $H^s(0,T) \subset \subset C^0([0,T])$
(or, more precisely, its vector-valued analogue),
holding for $s>1/2$.

\smallskip

We now claim that the above relations suffice to
take the limit $n\nearrow \infty$ 
in all equations of our system \eqref{incom}-\eqref{calore}.
To see this, we limit ourselves to consider the most troublesome
nonlinear terms, the other ones being in fact almost straighforward 
to treat. To start with, we note that
\begin{equation}\label{lim11}
  \nabla \fhi_n \otimes \nabla \fhi_n \to 
    \nabla \fhi \otimes \nabla \fhi 
   \quext{strongly in }\, C^0([0,T];L^1(\Omega)) 
\end{equation}
thanks to \eqref{cofhis}. Actually, proceeding as in
\eqref{co42}-\eqref{st42}, we can
check that $\dive (\nabla \fhi \otimes \nabla \fhi )$ lies
in $L^2(0,T;H)$, which allows 
for~\eqref{ns} to hold pointwise (almost everywhere).
Next, we notice that \eqref{covts} entails
$\vt_n \to \vt$ a.e.~in~$(0,T)\times \Omega$.
Such a property, combined with
\eqref{stkey2}, \eqref{stkey5}, and a generalized version of Lebesgue's
theorem, gives
\begin{equation}\label{lim12}
  \vt_n \to \vt \quext{strongly in }
   L^p(0,T;L^{(q+2)-}(\Omega)) \cap L^{(2q+2)-}(0,T;L^p(\Omega))
    \quext{for all }\,p\in[1,\infty).
\end{equation}
In particular, recalling \eqref{defiK}, we get
\begin{equation}\label{lim13}
  K(\vt_n) \to K(\vt) \quext{strongly in }
   L^p((0,T)\times \Omega)
    \quext{for a suitable }\,p>1.
\end{equation}
Hence, 
\begin{equation}\label{lim14}
  - \dive ( \kappa(\vt_n) \nabla \vt_n)
   = - \Delta K(\vt_n) 
   \to - \Delta K(\vt)
\end{equation}
at least in the sense of distributions. More precisely,
the limit function $K(\vt)$ lies in $\LDV$ in view of
\eqref{stkey5}, which allows the limit
of \eqref{calore} to hold as a relation in $L^2(0,T;V')$ as
specified by~\eqref{calore2}.

Moreover, recalling that $q\ge 2$ and using \eqref{cofhi},
\eqref{covus}, \eqref{cofhis} and \eqref{lim12},
we can easily check that
\begin{equation}\label{lim15}
  \vt_n ( \fhi_{n,t} + \vu_n \cdot \nabla \fhi_n )
   \to \vt ( \fhi_t + \vu \cdot \nabla \fhi )
   \quext{weakly in }
   L^p((0,T)\times \Omega)
    \quext{for a suitable }\, p>1.
\end{equation}
Next, thanks to \eqref{covus} and \eqref{comus}, 
we get (actually, something more is true)
\begin{equation}\label{lim16}
  | \nabla \vu_n |^2 + | \nabla \mu_n |^2
    \to | \nabla \vu |^2 + | \nabla \mu |^2
   \quext{strongly in }
   L^1((0,T)\times \Omega).
\end{equation}
Relations \eqref{lim11}-\eqref{lim14} permit us to let
$n\searrow \infty$ in all equations of the system
(with \eqref{calore} replaced by \eqref{calore2} 
in the limit, as noted above). In particular, the 
regularity properties \eqref{covus}-\eqref{covts} are 
a direct consequence of our argument. 
Finally, it is worth noting that the limit 
functions $\vu$, $\fhi$, $\vt$ also
satisfy the initial conditions \eqref{iniz}. Indeed,
they are continuous with respect to time with values in
suitable Banach spaces (and the corresponding
uniform estimates \eqref{covus}-\eqref{covts} hold for the
approximating sequences). Hence, we may
conclude that the limit quadruple $(\vu,\fhi,\mu,\vt)$
solves our system in the sense of Theorem~\ref{teo:main},
as desired.

%%%%%%%%%%%%%%%%%%%%%%%%%%%%%%%%%%%%%%%%%%%%%%%%%%%%%%%%%%%%%%%%%%%%%%%%%%%%%%%%%%%%%%%%%%%%%%%%%
%%%%%%%%%%%%%%%%%%%%%%%%%%%%%%%%%%%%%%%%%%%%%%%%%%%%%%%%%%%%%%%%%%%%%%%%%%%%%%%%%%%%%%%%%%%%%%%%%

%\section{Long-time behavior and attractors}
%\label{sec:long}

%%%%%%%%%%%%%%%%%%%%%%%%%%%%%%%%%%%%%%%%%%%%%%%%%%%%%%%%%%%%%%%%%%%%%%%%%%%%%%%%%%%%%%%%%%%%%%%%%
%%%%%%%%%%%%%%%%%%%%%%%%%%%%%%%%%%%%%%%%%%%%%%%%%%%%%%%%%%%%%%%%%%%%%%%%%%%%%%%%%%%%%%%%%%%%%%%%%

\section{Approximation and local existence}
\label{sec:appro}

In this section we give some highlights regarding a possible
approximation of system \eqref{incom}-\eqref{calore} and provide
a proof of local existence by means of a fixed point argument
of Schauder type. In order to reduce the length of the exposition
we leave most technical details to the reader, just limiting 
ourselves to outline the main steps of the procedure.

\smallskip

\noindent%
{\bf Regularized system.}~~%
For (small) $\epsi\in(0,1)$ we consider the following regularized statement:
\begin{align}\label{income}
  & \dive \vu = 0, \\
 \label{nse}
  & \vu_t + \vu \cdot \nabla \vu + \nabla p
    = \Delta \vu - \dive ( \nabla \fhi \otimes \nabla \fhi ),\\
 \label{CH1e}
  & \fhi_t + \vu \cdot \nabla \fhi = \Delta \mu, \\
 \label{CH2e}
  & \mu = -\epsi \Delta \fhi_t - \Delta \fhi + F\ee'(\fhi) - \vt, \\
 \label{caloree}
  & \vt_t + \vu\cdot \nabla \vt + \vt \big( \fhi_t + \vu \cdot \nabla \fhi \big) 
    - \dive(\kappa(\vt)\nabla \vt) = | \nabla \vu |^2 + T\ee\big( | \nabla \mu |^2 \big).
\end{align}  
The above system differs from the original one due in view of the 
additional term $-\epsi \Delta \fhi_t$ in \eqref{CH1e}, which provides
further parabolic regularity to $\fhi$, and of the 
{\sl truncation operator}\/ $T\ee$ in~\eqref{caloree}, where
\begin{equation}\label{trunc}
  T\ee(v):= \min\Big\{ \epsi^{-1}, \max\big\{ -\epsi^{-1}, v \big\} \Big\},
   \quext{for }\,v: (0,T)\times\Omega \to \RR,
\end{equation}  
which yields boundedness of the last term in the \rhs\ of \eqref{caloree}.
Moreover $F_\epsi$ is a smooth regularization of $F$ such that
$F_\epsi'$ is Lipschitz continuous. We assume that $F_\epsi$
still enjoys the coercivity property \eqref{hp:F1}.
We also truncate the initial temperature in such 
a way that
\begin{equation}\label{inize}
  %u\zee \in V, \qquad
  % \fhi\zee \in H^3(\Omega), \qquad
  \vt\zee \in H^1(\Omega)\cap L^\infty(\Omega), \qquad
   \epsi \le \vt\zee \le \epsi^{-1}~~\text{a.e.~in }\,\Omega.
\end{equation}  
Then, local existence for \eqref{income}-\eqref{caloree}, 
complemented with the initial data $\vu_0$, 
$\fhi_0$ and $\vt\zee$,
and with periodic boundary conditions, is proved via a 
fixed point argument detailed below. This is essentially
divided into three separate Lemmas. At first, 
we fix $\vt$ and $\vu$ in the Cahn-Hilliard system 
\eqref{CH1e}-\eqref{CH2e}.
\bele\label{lemma:fp1}
 Let $\epsi\in(0,1)$ and let $R>0$ be a number, depending on 
 the initial data and on $\epsi$, such that
 \begin{equation} \label{inizR}
   \| \varphi_0 \|_{H^3(\Omega)} 
    + \| \vu_0 \|_{V} 
    + \| \vt\zee \|_{H^1(\Omega)} 
    + \| \vt\zee \|_{L^\infty(\Omega)} 
   \le R.
 \end{equation}  
 Let also
 \begin{equation}\label{inizreg}
   \barvt \in \LDV, \quad
    \barvu \in L^2(0,T;L^4(\Omega)),
    \quext{with }\, \| \barvt \|_{\LDV} 
       + \| \barvu \|_{L^2(0,T;L^4(\Omega))} \le R.
 \end{equation}  
 Then there exist unique functions $\fhi$ 
 and $\mu$ satisfying, 
 a.e~in~$(0,T)\times\Omega$, the system
 \begin{align} \label{CH1fp}
   & \fhi_t + \barvu \cdot \nabla \fhi = \Delta \mu, \\
 \label{CH2fp}
   & \mu = -\epsi \Delta \fhi_t - \Delta \fhi + F\ee'(\fhi) - \barvt,
 \end{align}  
 with the initial condition $\fhi|_{t=0} = \fhi_0$.
 Moreover, the following regularity properties hold:
 \begin{align} \label{fhiee1}
    & \fhi\in H^1(0,T;H^3(\Omega)),\\
   \label{muee}
    & \mu \in L^2(0,T;H^2(\Omega)),\\
   \label{fhiee2}
    & \| \fhi \|_{H^1(0,T;H^3(\Omega))} 
     + \| \mu \|_{L^2(0,T;H^2(\Omega))} \le Q_1(R,T).
 \end{align}  
 Here and below, $Q_i:(\RR^+)^2\to \RR^+$, $i=1,2,\dots$,
 are computable functions, increasingly monotone 
 in each of their arguments, whose expression may
 additionally depend on $\epsi$.
\enle
\begin{proof}
We just give the highlights. Actually, once $\barvu$ and $\barvt$ are 
assigned, \eqref{CH1fp}-\eqref{CH2fp} is a semilinear pseudo-parabolic
system with Lipschitz nonlinearity. Hence, existence is standard. 
For example, it could be proved by relying on a Faedo-Galerkin
scheme, or on a time-discretization argument. The 
a-priori estimates corresponding to the regularity conditions \eqref{fhiee1}-\eqref{muee}
are the following ones: first, one reproduces the energy estimate for the complete
system by testing \eqref{CH1fp} by $\mu$ and \eqref{CH2fp} by $\fhi_t$. Noticing
that 
\begin{equation}\label{stfp11}
  \io \mu \, \barvu\cdot \nabla \fhi
   = - \io \fhi \, \barvu \cdot \nabla \mu
   \le \| \fhi \|_{L^4(\Omega)} \| \barvu \|_{L^4(\Omega)} \| \nabla \mu \|
   \le \frac12 \| \nabla \mu \|^2 + c \| \fhi \|_{V}^2 \| \barvu \|_{L^4(\Omega)}^2,
\end{equation}
an estimate follows by using \eqref{inizreg} and applying the Gronwall Lemma.
With this estimate at disposal we can test \eqref{CH2fp} by 
$\Delta^2 \fhi_t$. Thanks to the fact that $\mu,\barvt\in \LDV$ due
to the energy estimate and to \eqref{inizreg}, and using 
the Lipschitz continuity of $F_\epsi'$, it is then not difficult
to obtain~\eqref{fhiee1}. 

Subsequently, noting that the \lhs\ of \eqref{CH1fp}
lies in $L^2(0,T;H)$, 
by elliptic regularity we obtain \eqref{muee}. Then,
relation \eqref{fhiee2} is also a direct consequence of the a-priori estimates,
as one can see by writing them in a quantitative way. Finally, uniqueness 
in the regularity class specified by \eqref{fhiee1}-\eqref{muee} can be proved
by a standard contractive argument. Namely, one may test (the
difference of) \eqref{CH1fp} by $\fhi_1-\fhi_2$ (where $(\fhi_1,\mu_1)$
and $(\fhi_2,\mu_2)$ are two solutions) and (the difference of)
\eqref{CH2fp} by $\Delta(\fhi_1-\fhi_2)$ and perform standard
calculations.
\end{proof}
\bele\label{lemma:fp2}
 Let us assume that the hypotheses of~{\rm Lemma~\ref{lemma:fp1}} are
 satisfied, and let $\fhi$, $\mu$ be the functions provided 
 by~{\rm Lemma~\ref{lemma:fp1}}. Then, there exists 
 a unique function $\vu$ such that
 \begin{align} \label{uee1}
    & \vu\in H^1(0,T;H) \cap \LIV \cap L^2(0,T;H^2(\Omega)), \\
   \label{uee2}
    & \| \vu \|_{H^1(0,T;H)} + \| \vu \|_{\LIV}
     + \| \vu \|_{L^2(0,T;H^2(\Omega))}
    \le Q_2(R,T).
 \end{align}  
 Moreover, $\vu$ satisfies, 
 a.e~in~$(0,T)\times\Omega$, the system
 \begin{align}\label{incomfp}
   & \dive \vu = 0, \\
  \label{nsfp}
   & \vu_t + \vu \cdot \nabla \vu + \nabla p - \Delta \vu
    = - \dive ( \nabla \fhi \otimes \nabla \fhi ),
 \end{align}  
 with the initial condition $\vu|_{t=0} = \vu_0$.
\enle
\begin{proof}
Also in this case we just give the highlights. Actually, as 
$\fhi$ is given satisfying \eqref{fhiee1} and \eqref{fhiee2},
it is clear that the \rhs\ of \eqref{nsfp} 
lies in $L^2(0,T;H)$. Hence, existence
and uniqueness of a solution satisfying \eqref{uee1}
follow from the general theory of Navier-Stokes systems 
(cf., e.g., \cite{rob} or \cite{Tem}).
Moreover, writing explicitly the a-priori bounds, 
one also immediately gets \eqref{uee2}, 
where, in principle, the expression of
$Q_2$ may also depend on a suitable norm of $\fhi$.
However, thanks to \eqref{fhiee2}, $Q_2$ can in fact
be written as a (computable and monotone) function of $R$ and $T$.
\end{proof}
\noindent%
Finally, we come to the ``heat'' equation, which is 
a little bit more involved to deal with:
\bele\label{lemma:fp3}
 Let the assumptions of~{\rm Lemma~\ref{lemma:fp1}} hold
 and let $\fhi$, $\mu$, $\vu$ be the functions provided 
 by~{\rm Lemmas~\ref{lemma:fp1}, \ref{lemma:fp2}}. Then, there exists 
 a unique function $\vt$ such that
 \begin{align} \label{vtee1}
    & \vt\in H^1(0,T;H) \cap \LIV \cap L^2(0,T;H^2(\Omega)) 
             \cap L^\infty((0,T)\times \Omega), \qquad
     \vt > 0~~\text{a.e.~in~$(0,T)\times\Omega$,} \\
   \label{vtee2}
    & \| \vt \|_{H^1(0,T;H)} 
     + \| \vt \|_{\LIV}
     + \| \vt \|_{L^2(0,T;H^2(\Omega))}
    \le Q_3(R,T).
 \end{align}  
 Moreover, $\vt$ satisfies, 
 a.e~in~$(0,T)\times\Omega$, 
 \begin{equation}\label{calorefp} 
   \vt_t + \vu\cdot \nabla \vt + \vt \big( \fhi_t + \vu \cdot \nabla \fhi \big) 
    - \dive(\kappa(\vt)\nabla \vt) = | \nabla \vu |^2 + T\ee\big( | \nabla \mu |^2 \big),
 \end{equation}  
 with the initial condition $\vt|_{t=0} = \vt\zee$.
\enle
\begin{proof}
Equation \eqref{calorefp} enjoys the quasilinear structure
\begin{equation}\label{calorefp2}
  \vt_t + \vu\cdot \nabla \vt + m_1 \vt 
   - \Delta K(\vt) = f,
\end{equation}  
where $K$ was defined in \eqref{defiK}, and,
in view of estimates \eqref{uee2}, \eqref{fhiee2}
and of interpolation, it is not difficult to infer
\begin{align} \label{regfp1}
  & f := | \nabla \vu |^2 + T\ee\big( | \nabla \mu |^2 \big)
   \in L^p((0,T)\times\Omega) 
   \quext{for all }\ p\in (1,\infty),\\
 \label{regfp2}
  & m_1 := \fhi_t + \vu \cdot \nabla \fhi
   \in L^1(0,T;L^\infty(\Omega)) \cap L^2(0,T;L^p(\Omega))
   \quext{for all }\ p\in (1,\infty).
\end{align}  
Hence, existence of solutions to
the initial-boundary value problem for \eqref{calorefp2}
follows, as before, from standard techniques (as
before, one could use time discretization or
Faedo-Galerkin approximation).

Moreover, testing \eqref{calorefp2}
by $\vt$ and performing simple calculations, we
obtain the a-priori estimates leading to the 
regularity $\vt\in L^\infty(0,T;H) \cap L^2(0,T;V)$. 
By virtue of the high summability of $f$ and $m_1$,
a standard application of Moser's iteration technique
(see, e.g., \cite[Chap.~5]{LSU})
yields $\vt\in L^\infty((0,T)\times\Omega)$. 

Strict positivity of $\vt$ follows from 
the maximum principle. Next, to deduce 
\eqref{vtee1}-\eqref{vtee2} one tests \eqref{calorefp}
by $K(\vt)_t$. Recalling \eqref{hp:kappa}, and using the $L^\infty$-bound
for $\vt$ coming from Moser's iterations,
we then get the estimate
\begin{align}\no
  & \io ( 1 + \vt^q ) | \vt_t |^2
   + \frac12 \ddt \| \nabla K(\vt) \|^2
   = \io (f - \vu\cdot \nabla \vt - m_1 \vt ) ( 1 + \vt^q ) \vt_t\\
 \no 
  & \mbox{}~~~~~ 
  \le \frac12 \| \vt_t \|^2
   + c \big\| f - \vu \cdot \nabla \vt - m_1 \vt \big\|^2 
    \big( 1 + \| \vt \|_{L^\infty(\Omega)}^{2q} \big)\\
 \label{stfp21}
  & \mbox{}~~~~~ 
  \le \frac12 \| \vt_t \|^2
   + c \| f \|^2 + c \big( \| \vu \|_{L^\infty(\Omega)}^2
   + \| m_1 \|_{L^4(\Omega)}^2 \big) \| \vt \|_V^2.
\end{align}  
Hence, noting that $\| \vt \|_V \le c \| K(\vt) \|_V$
by \eqref{hp:kappa}, conditions $\vt \in \HUH$
and $K(\vt) \in \LIV$ follow from Gronwall's lemma. 
Consequently we also have $\vt \in \LIV$.
In addition to that, viewing \eqref{calorefp2} as a 
time-dependent family of elliptic problems 
(for the variable $K(\vt)$) with $L^2$-data
and applying standard regularity results,
we obtain that $K(\vt) \in \LDHD$.

Let now $k:=K(\vt)$ and denote
as $\eta_q$ the inverse function of $K$
over $[0,+\infty)$. Observe that 
$\eta_q'$ and $\eta_q''$ are uniformly bounded. Hence
\begin{equation}\label{Kk}
  \Delta \vt 
   = \Delta \eta_q(k)
   = \eta_q'(k) \Delta k + \eta_q''(k) | \nabla k |^2
\end{equation}  
belongs to $\LDH$ thanks to the fact that 
$k\in\LDHD \cap \LIV$. This entails \eqref{vtee1}.
As before, as one writes explicitly the estimates 
leading to \eqref{vtee1}, also \eqref{vtee2} follows.

Finally, to show uniqueness, we take two 
solutions $\vt_1$ and $\vt_2$ to \eqref{calorefp}
with the same initial datum and the same $\fhi,\mu,\vu$.
Then, the difference $\tilde\vt:=\vt_1-\vt_2$ solves
\begin{equation}\label{calorefp2d}
  \tilde\vt_t + \vu\cdot \nabla \tilde\vt + m_1 \tilde\vt 
   - \Delta( K(\vt_1) - K(\vt_2) ) = 0.
\end{equation}  
Then, testing the above by $\sign \tilde\vt$ (more precisely,
one should first take an approximation of the sign function
and then pass to the limit), using monotonicity
of $K$, the incompressibility \eqref{incomfp}, and the
periodic boundary conditions, we infer
\begin{equation}\label{stfp22}
  \ddt \| \tilde\vt \|_{L^1(\Omega)}
   \le \| m_1 \|_{L^\infty(\Omega)} \| \tilde\vt \|_{L^1(\Omega)},
\end{equation}  
whence uniqueness follows from Gronwall's Lemma recalling
the first condition in~\eqref{regfp2}. 
\end{proof}
\noindent%
With the three lemmas at disposal, we can make explicit
our fixed-point argument
\bete\label{teo:schauder}
 Let $\epsi\in(0,1)$ and let us assume \eqref{trunc} and \eqref{inize}.
 Then there exist a time $T_0$ (depending on $\epsi$ 
 and on the initial data) and at least one quadruple 
 $(\vu,\fhi,\mu,\vt)$ such that
 \begin{align}\label{regouee}
   & \vu \in H^1(0,T_0;H) \cap L^\infty(0,T_0;V)
     \cap L^2(0,T_0;H^2(\Omega)), \\
  \label{regofhiee}
   & \fhi \in H^1(0,T_0;H^3(\Omega)),\\
  \label{regomuee}
   & \mu \in L^2(0,T_0;H^2(\Omega)),\\  
  \label{regovtee}
   & \vt\in H^1(0,T_0;H)  \cap L^\infty(0,T_0;V)
      \cap L^2(0,T_0;H^2(\Omega)) \cap L^\infty((0,T_0)\times \Omega), 
 \end{align}  
 with $\vt > 0$ a.e.~in~$(0,T_0)\times \Omega$,
 satisfying system~{\rm \eqref{income}-\eqref{caloree}}
 a.e.~in~$(0,T_0)\times\Omega$ 
 and complying with the initial conditions
 \begin{equation}\label{cauchyee}
   \vu|_{t=0} = \vu_0, \qquad
   \fhi|_{t=0} = \fhi_0, \qquad
   \vt|_{t=0} = \vt\zee.
 \end{equation}  
\ente
\begin{proof}
Given $\epsi>0$, we truncate the initial temperature
as specified in \eqref{inize}. Then, we choose $R>0$ correspondingly
(cf.~\eqref{inizR}). Hence, we can consider the closed ball
(cf.~\eqref{inizreg})
\begin{equation}\label{ball}
  \calB:= \big\{ (\barvt,\baru): \| \barvt \|_{L^2(0,T_0;V)}
       + \| \barvu \|_{L^2(0,T_0;L^4(\Omega))} \le R \big\},
\end{equation}  
where $T_0\in(0,T]$ will be chosen later on. Notice that the chosen
radius $R$ depends only on the initial data (and on $\epsi$
by the truncation applied to $\vt_0$). Let us consider
the fixed point map (also depending on $\epsi$, of course)
\begin{equation}\label{fixmap}
  \calT: \calB \to L^2(0,T_0;V) \times L^2(0,T_0;L^4(\Omega)), \qquad
   \calT:(\barvt,\baru)\mapsto (\vt,\vu).
\end{equation}  
We aim at applying the Schauder fixed point theorem to
the above map, for a sufficiently small choice of the 
final time $T_0>0$. To this aim, we can observe the following:
\begin{itemize}
\item[(a)]~~The map~$\calT$ is continuous: this follows from
 the fact that the fixed point equations \eqref{CH1fp}-\eqref{CH2fp},
 \eqref{nsfp} and \eqref{calorefp} only contain Lipschitz 
 or locally Lipschitz nonlinearities. To give a formal proof
 (which is omitted for brevity), one could just
 put together (and refine a bit) the contractive arguments 
 used to prove uniqueness in the three fixed point Lemmas.
\item[(b)]~~The map~$\calT$ is compact: this follows immediately
 from \eqref{vtee1}, \eqref{uee1} and the Aubin-Lions lemma.
 Indeed, the output of the map $\calT$ lies in a bounded set
 of a space which is compactly embedded into 
 $L^2(0,T_0;V) \times  L^2(0,T_0;L^4(\Omega))$.
\item[(c)]~~The map~$\calT$ takes values into $\calB$. 
 Indeed, thanks to \eqref{vtee2}, \eqref{uee2}, 
 and the continuous embedding $V\subset L^4(\Omega)$,
 we get
 \begin{equation}\label{uvtfp}
   \| \vu \|_{L^\infty(0,T_0;L^4(\Omega))} 
    + \| \vt \|_{L^\infty(0,T_0;V)}
    \le Q_4(R,T_0),
 \end{equation}  
 for some function $Q_4(R,T_0)$, whence
 \begin{equation}\label{uvtfp2}
   \| \vu \|_{L^2(0,T_0;L^4(\Omega))} 
    + \| \vt \|_{L^2(0,T_0;V)}
    \le T_0^{1/2} Q_4(R,T_0) \le R,
 \end{equation}  
 provided $T_0$ is small enough. 
\end{itemize}
In view of the above conditions (a)-(c) the assumptions of
the Schauder fixed point argument are satisfied,
whence (at least) one solution to \eqref{income}-\eqref{caloree}
exists. The regularity conditions \eqref{regouee}-\eqref{regovtee}
follow immediately from the above three Lemmas.
Theorem~\ref{teo:schauder} is proved.
\end{proof}
\noindent%
We conclude this section with three observations aimed
at clarifying why the present approximation-fixed point 
argument is compatible with the a-priori estimates of
the previous section.
\beos\label{remfp1}
 As usual, the solution provided by the fixed point argument is 
 local in time and the final time $T_0$ may depend on $\epsi$ 
 and be smaller as smaller is $\epsi$. However, the a-priori
 estimates performed in the previous section are uniform
 with respect to time. Hence, thanks to standard 
 extension arguments, the (approximate) solution 
 turns out to be defined, in fact, on the whole 
 reference interval~$(0,T)$,
 and the same will hold for the limit solution.
\eddos
\beos\label{remfp2}
 It is also worth noting that neither the regularizing term
 $-\epsi\Delta\fhi_t$ added on the \rhs\ of \eqref{CH2e} nor 
 the truncation operator on the \rhs\ of \eqref{caloree} 
 really interfer with the a-priori estimates 
 of the previous section, which turn out to be independent
 of the approximation parameter $\epsi$. Actually,
 $-\epsi\Delta\fhi_t$ just gives some more information,
 vanishing as $\epsi$ goes to $0$, in the energy and subsequent
 bounds. On the other hand, the truncation operator 
 yields some (positive) remainder term on the \lhs\ of the 
 energy estimate, coming from the fact that, as one tests \eqref{CH1e} 
 by $\mu$ and \eqref{caloree} by $1$, the contribution of $\mu$ does 
 not vanish completely. It is immediate to see that 
 this additional term compensates exactly the lack of information
 one gets in the subsequent temperature estimate 
 (due to the fact that now only the truncated 
 $|\nabla \mu|^2$ appears on the \lhs\ of \eqref{co31}).
 Hence, \eqref{st32} can still be obtained.
\eddos
\beos\label{remfp3}
 Finally, we notice that the regularity class 
 \eqref{regouee}-\eqref{regovtee} is not really sufficient
 in order for the estimates of the previous section to 
 be rigorous. For instance, in principle we have no information
 on the term $\mu_t$, which would be needed as we perform the 
 ``Key estimate''. However, it is easy to realize that
 the regularity given for the local approximate solution
 is just the outcome of the fixed point argument and, hence,
 is not at all optimal. Further regularity properties 
 can be standardly proved by working separately
 on the equations of the approximate system and 
 performing simple bootstrap arguments.
 This procedure may involve some boring technical details
 (and may also require some further regularization of the 
 initial data) and is omitted for brevity.
 However, it is worth noting at least 
 that, for what concerns the ``heat'' equation (which contains
 the more delicate nonlinearities), the regularity obtained
 in the fixed-point argument is sufficient.
 So, it would be sufficient to bootstrap regularity 
 for $\vu$ and $\fhi$.
\eddos

 %%%%%%%%%%%%%%%%%%%%%%%%%%%%%%%%%%%%%%%%%%%%%%%%%%%%%%%%%%%%%%%%%%%%%%%%%%%%%%%%%%%%%%%%%%%5

 \appendix

\section{Appendix}
 \label{sec:appe}

We prove here the two-dimensional interpolation-embedding
inequality used in the proof of existence. It is worth noting
that the result is independent of the use of 
periodic boundary conditions
and holds in any smooth bounded subset of $\RR^2$.
\bele\label{lemma:yudo}
 Let $\calO$ a smooth bounded domain in $\RR^2$. 
 Then, there exists $c>0$ depending only on~$\calO$ 
 such that
 \begin{equation}\label{yudo}
   \| \xi \|_{H^1(\calO)'} \le c 
     \Big( 1 + \| \xi \|_{L^1(\calO)} 
          \log^{1/2} \big(e + \| \xi \|_{L^2(\calO)} \big) \Big)
 \end{equation}
 for any $\xi\in L^2(\calO)$.
\enle
\begin{proof}
We start by recalling (see, e.g., \cite[(17), p.~479]{Trud})
that, for all $p\in [1,\infty)$,
\begin{equation}\label{best}
  \| v \|_{L^p(\calO)} 
    \le c p^{1/2} \| v \|_{H^1(\calO)} 
   \quext{for all }\, v\in H^1(\calO).
\end{equation}
where the constant $c>0$ can be taken {\sl independent}\/
of $p$. The above inequality makes precise the rate of explosion
of the embedding constant of $H^1$ into $L^p$ as $p$ becomes
large. As before, the value of $c$ can vary in the computations
below; in any case, $c$ will always be intended to be
independent of~$p$.

Let us identify, as before, $L^2(\calO)$ with its dual in
such a way that $L^2(\calO)$ can be (compactly and continuously)
embedded into $H^1(\calO)'$. Then, given $\eta \in L^2(\calO)$, 
we have, for $c>0$ as above,
\begin{equation}\label{yudo11}
  \frac{\duav{ \eta, v}}{\| v \|_{H^1(\calO)}}
   \le c p^{1/2} \frac{\duav{ \eta, v}}{\| v \|_{L^p(\calO)}}
   \le c p^{1/2} \| \eta \|_{L^{p'}(\calO)},
\end{equation}
for all nonzero $v\in H^1(\calO)$ and where $p'$ is the 
conjugate exponent to $p$. Hence,
\begin{equation}\label{yudo12}
  \| \eta \|_{H^1(\calO)'}
     \le c p^{1/2} \| \eta \|_{L^{p'}(\calO)}
     = c p^{1/2} \| \eta \|_{L^{\frac{p}{p-1}}(\calO)}.
\end{equation}
Taking $p\ge 2$ and using interpolation, it then follows
\begin{equation}\label{yudo13}
  \| \eta \|_{H^1(\calO)'}
     \le c p^{1/2} 
     \| \eta \|_{L^{1}(\calO)}^{\frac{p-2}p}
          \| \eta \|_{L^{2}(\calO)}^{\frac{2}p}.
\end{equation}
Let us temporarily assume that $ \| \eta \|_{L^{2}(\calO)} = 1$. 
Then, squaring, we obtain
\begin{equation}\label{yudo14}
  \| \eta \|_{H^1(\calO)'}^2
     \le c p
     \| \eta \|_{L^{1}(\calO)}^{\frac{2(p-2)}p}
     = c p \| \eta \|_{L^{1}(\calO)}^{2}
          \| \eta \|_{L^{1}(\calO)}^{-\frac{4}p}.
\end{equation}
Now we use the so-called Yudovich' trick 
(see, e.g., \cite{Yud}), 
namely we optimize the above \rhs\
with respect to $p\in [2,\infty)$. To this aim,
let us consider the function
\begin{equation}\label{fp}
  f:[2,\infty)\to (0,\infty),
   \qquad f(p) := p A^{-\frac{4}p},
\end{equation}
where $A>0$ is given. Then, clearly,
\begin{equation}\label{fpprimo}
  f'(p) = A^{-\frac{4}p} \Big( 1 + \frac4p \log A \Big).
\end{equation}
Now, we have to distinguish between two cases. 
Firstly, if $A \ge e^{-1/2}$, then $f$ is increasing
over $[2,+\infty)$, whence its minimum is achieved 
for $p=2$:
\begin{equation}\label{yudo21}
  \min f = f(2) = 2 A^{-2} \le 2 e.
\end{equation}
On the other hand, if $A\in (0,e^{-1/2})$, then 
$-4 \log A$ (the zero of $f'$) is strictly 
larger than~$2$, whence
\begin{equation}\label{yudo22}
  \min f = f( - 4 \log A ) = 
   - 4 A^{\frac1{\log A}} \log A  
   = - 4 e \log A.
\end{equation}
Let us now choose $A = \| \eta \|_{L^1(\calO)}$ for 
$\eta \in L^2(\calO)$ with $\| \eta \|_{L^2(\calO)}=1$.
Then,
\begin{equation}\label{yudo23}
  \| \eta \|_{H^1(\calO)'}^2
   \le \begin{cases}
           2 e c \| \eta \|_{L^1(\calO)}^2 
	     & \text{~~if }\, \| \eta \|_{L^1(\calO)} \ge e^{-1/2}\\
	   - 4 e c \| \eta \|_{L^1(\calO)}^2 
             \log \big( \| \eta \|_{L^1(\calO)} \big)  
	     & \text{~~if }\, \| \eta \|_{L^1(\calO)} < e^{-1/2}.
	\end{cases}    
\end{equation}
Let us now take any nonzero $\xi \in L^2(\calO)$ and apply the above
to $\eta = \xi / \| \xi \|_{L^2(\calO)}$. If 
\begin{equation}\label{yudo24}
  \| \eta \|_{L^1(\calO)} \ge e^{-1/2}, 
   \quext{i.e.~}\, \| \xi \|_{L^2(\calO)} \le e^{1/2} \| \xi \|_{L^1(\calO)},
\end{equation}
then it follows from the first \eqref{yudo23} that
\begin{equation}\label{yudo1}
  \| \xi \|_{H^1(\calO)'}^2
     \le 2 e c \| \xi \|_{L^1(\calO)}^2, 
\end{equation}
and, in particular, \eqref{yudo} holds. On the other hand, if 
\begin{equation}\label{yudo25}
  \| \eta \|_{L^1(\calO)} < e^{-1/2}, 
   \quext{i.e.~}\, \| \xi \|_{L^2(\calO)} > e^{1/2} \| \xi \|_{L^1(\calO)},
\end{equation}
then from the second \eqref{yudo23} we obtain
\begin{equation}\label{yudo26}
  \| \xi \|_{H^1(\calO)'}^2 
     \le 4 e c \| \xi \|_{L^1(\calO)}^2 
     \big( \log \| \xi \|_{L^2(\calO)} - \log \| \xi \|_{L^1(\calO)} \big).
\end{equation}
Here, we have to distinguish again some cases. First,
if it is both $\| \xi \|_{L^2(\calO)} \ge 1$ 
and $\| \xi \|_{L^1(\calO)} \ge 1$, then 
\eqref{yudo26} is continuated as
\begin{equation}\label{yudo26b}
  \| \xi \|_{H^1(\calO)'}^2 
     \le 4 e c \| \xi \|_{L^1(\calO)}^2 
      \log \| \xi \|_{L^2(\calO)}, 
\end{equation}
whence \eqref{yudo} follows. Second, if 
$\| \xi \|_{L^2(\calO)} \ge 1$ 
and $\| \xi \|_{L^1(\calO)} < 1$, then, observing 
that $ 0 \le - r \log r \le c$ for all $r\in(0,1)$, we
get
\begin{equation}\label{yudo26c}
  \| \xi \|_{H^1(\calO)'}^2 
     \le 4 e c \| \xi \|_{L^1(\calO)}^2 \log \| \xi \|_{L^2(\calO)}
       - 4 e c \| \xi \|_{L^1(\calO)}^2 \log \| \xi \|_{L^1(\calO)} 
     \le 4 e c \| \xi \|_{L^1(\calO)}^2 \log \| \xi \|_{L^2(\calO)} + c,
\end{equation}
and we still have \eqref{yudo}. Finally, if it is both
$\| \xi \|_{L^2(\calO)} < 1$ and $\| \xi \|_{L^1(\calO)} < 1$
then, we simply observe that 
$\| \xi \|_{H^1(\calO)'}^2 \le c \| \xi \|_{L^2(\calO)}^2 \le c$,
which concludes the proof.
\end{proof}
%

%%%%%%%%%%%%%%%%%%%%%%%%%%%%%%%%%%%%%%%%%%%%%%%%%%%%%%%%%%%%%%%%%%%%%%%%%%%%%%%%%%%%%%%%%%%%5

\bibliographystyle{amsplain}

%%%%%%%%%%%%%%%%%%%%%%%%%%%%%%%%%%%%%%%%%%%%%%%%%%%%%%%%%%%%%%%%%%%%%%%%%%%%%%%%%%%%%%%%%%%%%

\end{document}